\documentclass[onefignum,onetabnum,noabbrev]{siamonline171218}

\usepackage{multirow}
\usepackage{siunitx}

\newcommand{\vicadd}[1]{{#1}}
\newcommand{\vicdelete}[1]{}
\newcommand{\vicreplace}[2]{{#2}}

\usepackage{amsfonts}
\usepackage{amssymb}
\usepackage{bm}
\usepackage{todonotes}
\usepackage{microtype}
\usepackage{graphicx}
\usepackage{tikz}
\usepackage{epstopdf}
\usepackage{algpseudocodex}
\usepackage{enumitem}
\usepackage{subcaption}
\usepackage{bbm}
\usepackage{booktabs}
\setlist[enumerate]{leftmargin=.5in}
\setlist[itemize]{leftmargin=.5in}
\ifpdf
  \DeclareGraphicsExtensions{.eps,.pdf,.png,.jpg}
\else
  \DeclareGraphicsExtensions{.eps}
\fi

\crefformat{equation}{Eq.~(#2#1#3)}

\DeclareMathOperator*{\argmax}{arg\,max}

\newcommand{\fb}{\mathbf{f}}\newcommand{\Fb}{\mathbf{F}}\newcommand{\yb}{\mathbf{y}}
\newcommand{\gb}{\mathbf{g}}\newcommand{\Ex}{\mathbb{E}}\newcommand{\Prob}{\mathbb{P}}
\newcommand{\Xspace}{\mathcal{X}}
\newcommand{\Uspace}{\mathcal{U}}
\newcommand{\ParFront}{\mathcal{P}^*}
\newcommand{\ParSet}{\mathcal{P}^*_\Xspace}
\newcommand{\mb}{\mathbf{m}}

\renewcommand{\epsilon}{\varepsilon}

\newcommand{\HV}{\mathop{\mathrm{HV}}}

\headers{MOOUU with Conditional Pareto Fronts}{V. Trappler, C. Helbert, R. Le Riche}

\title{Multiobjective Optimization under Uncertainties using Conditional Pareto Fronts}

\author{Victor Trappler\thanks{Centrale Lyon, CNRS UMR 5208, Institut Camille Jordan. Now with Mines Saint-Etienne, Univ Clermont Auvergne, INP Clermont Auvergne, CNRS, UMR 6158 LIMOS, F - 42023 Saint-Etienne, France}
\and Céline Helbert\thanks{Centrale Lyon, CNRS UMR 5208, Institut Camille Jordan, 36
Avenue Guy de Collongue, 69134 Écully, France} \and Rodolphe Le Riche\thanks{
Laboratoire d’Informatique, de Modélisation et d’Optimisation des Systèmes (LIMOS), 63178 Aubiere, France}}
\begin{document}
\maketitle
\begin{abstract}
In this work, we propose a novel method to tackle the problem of multiobjective
optimization under parametric uncertainties, by considering the Conditional
Pareto Sets and Conditional Pareto Fronts. Based on those quantities, we can
define the probability of coverage of the Conditional Pareto Set which can be
interpreted as the probability for a design to be optimal in the Pareto sense.
Due to the computational cost of such an approach, we introduce an Active
Learning method based on Gaussian Process Regression in order to improve the
estimation of this probability, which relies on a reformulation of the \vicreplace{EHVI}{Expected Hypervolume Improvement (EHVI)}. We
illustrate those methods on a few toy problems of moderate dimension, and on the
problem of designing a cabin to highlight the differences in solutions brought
by different formulations of the problem.
\end{abstract}



\section{Introduction}
In many industrial or scientific studies, the task of finding suitable
parameters is often formulated as an optimization problem where the controlled
parameters should be the best according to some specific criterion, which
represents the cost of making such a decision. When several objectives are
considered simultaneously, the decision maker has different possibilities to
tackle the problem. For instance, one might want to combine the objectives
through a convex combination, or to optimize one of the objective while adding
inequality constraints on the other ones. A more general approach is to look for
all the best possible \vicreplace{compromises}{tradeoffs} between those objectives, namely the Pareto
front and its pre-image the Pareto set, which are the multiobjective
counterparts of the optimum and optimizers, respectively.

Most popular global multiobjective optimization methods, either in the mono or
multiobjective case require a large number of evaluations of the objective
function. This is due both to the stochastic sampling that is needed to escape
from local optima, and to the cost of populating the Pareto set. Examples of
such methods are the multiobjective version of
CMA-ES~\cite{igel_covariance_2007,toure2019uncrowded},
NSGA-II~\cite{deb_fast_2002}, or any scalarizations of the multiobjective
problem that translates into a series of mono-objective problems
\cite{zhang2020random}. However, in many practical cases, a single evaluation of
the cost function is expensive because it may require expensive computer
simulations or physical experiments. In that case, when trying to solve an
optimization problem, being able to limit the total number of function
evaluations is critical. This motivates the use of Bayesian Optimization (see
\cite{frazier2018tutorial,garnett_bayesian_2023,wang_recent_2023} for extensive
reviews), where we assume specific priors on the unknown objective functions,
and add points sequentially to the design of experiment according to some
progress measure called the acquisition function. In the past decades, this
derivative-free method has been applied to various problems where the budget of
simulation is limited and/or costly, in industrial contexts for instance \cite{gaudrie_modeling_2020}, hyperparameters tuning in Machine
Learning \cite{klein_fast_2017}, drug discovery \cite{colliandre_bayesian_2024}{or quantitative microbiological risk assessment \cite{basak_multipathogen_2024}.

One specific case of interest is when the objective function is not only a
function of the control variables, but also depends on environmental variables
which represents some uncertainties in the model. This separation of control
variables and uncertain variables has been used mostly when the underlying
objective function is deterministic, and the choice of both variables is up to
the user~\cite{lehman_designing_2004,trappler_robust_2021,elamri_sampling_2023}.

In this work, we propose to tackle the problem of multiobjective optimization
under uncertainties through the notion of conditional Pareto front. Since the
Pareto front and Pareto set are random quantities that depend solely on the
environmental variables, we can consider the probability of coverage of the
Pareto set as a robustness measure. This probability can be expensive to
calculate, and this limitation can become even more critical when done within an
optimization procedure. That is why we adapt to the multiobjective case the
acquisition function introduced in \cite{ginsbourger_bayesian_2014} in order to
improve the estimation of the Conditional Pareto Fronts using Gaussian Process
Regression.

\subsection{Multiobjective Optimization}
\label{ssec:moo}
Let $\Xspace \subset \mathbb{R}^{n_x}$\vicreplace{, and let the objective function}{. The objective function is defined as}
\begin{equation}
\begin{array}{rcl}
\fb: \Xspace & \longrightarrow& \mathbb{R}^d \\
x & \longmapsto &\fb(x) = (f_1(x),\dots,f_d(x))
\end{array}\,.
\end{equation}

The function $\fb$ maps $x\in \Xspace$ to a real vector of dimension $d$. Each
of the components of $\fb(x)$ represents an objective that we wish to minimize.
In practical cases, $d$ is rarely larger than 3, as the case $d\geq 4$ calls for
specific methods of Many Objective Optimization
\cite{fleming_manyobjective_2005,ishibuchi_evolutionary_2008,cai_kernelbased_2022}.
We are interested formally in the following multiobjective optimization problem
\begin{align}
    \min_{x \in \Xspace} \fb(x) &= (f_1(x),\dots,f_d(x))\,.
    \label{eq:moo}
\end{align}

In Multiobjective Optimization (MOO), the decision maker wants to choose a
 design $x\in \Xspace$ such that all the components of the objective function $\fb$ are
 minimized. Assuming that $\fb$ is continuous and $\Xspace$ compact, we call the
 ideal point the most optimistic objective specification
 $\mathbf{y}_{\mathrm{ideal}} = \left(\min f_1,\dots,\min f_d\right)$. On the
 other hand, the most pessimistic objective is usually called the nadir:
 $\mathbf{y}_{\mathrm{nadir}} = \left(\max f_1,\dots,\max f_d\right)$. In most
 problems the objectives are competing and thus do not share the same minimizer:
 the ideal cannot be reached. 

In order to compare the performances of different control vectors, we introduce
a partial order called the weak dominance order in $\mathbb{R}^d$. For $\yb =
(y_1,\dots,y_d)$ and $\yb' = (y'_1,\dots,y'_d)$, we say that $\yb$ dominates
$\yb'$ and we note $\yb \prec \yb'$ if $\forall i,\, y_i \leq y_i' \text{ and }
\exists j $ such that $y_j <y_j'$. In other words, $\yb \prec \yb'$ means that
$\yb$ is at least as good as $\yb$ in all objectives, and strictly better in at
least one. One evaluation $\fb(x)$ of the function partitions the objective
space into dominated, non-dominated and incomparable regions, as shown in
\cref{fig:domination_relation}. 

\begin{figure}
    \centering
    \scalebox{.8}{\tikzset{every picture/.style={line width=0.75pt}} 

\begin{tikzpicture}[x=0.75pt,y=0.75pt,yscale=-1,xscale=1]

\draw [color={rgb, 255:red, 208; green, 2; blue, 27 }  ,draw opacity=1 ]   (318.65,105) -- (432.05,105) ;
\draw  (191.8,179.04) -- (446.6,179.04)(217.28,10.2) -- (217.28,197.8) (439.6,174.04) -- (446.6,179.04) -- (439.6,184.04) (212.28,17.2) -- (217.28,10.2) -- (222.28,17.2)  ;
\draw  [draw opacity=0][fill={rgb, 255:red, 221; green, 99; blue, 99 }  ,fill opacity=0.58 ] (318,20.2) -- (432.05,20.2) -- (432.05,105) -- (318,105) -- cycle ;
\draw  [draw opacity=0][fill={rgb, 255:red, 126; green, 211; blue, 33 }  ,fill opacity=0.51 ] (201.8,105) -- (318,105) -- (318,191.4) -- (201.8,191.4) -- cycle ;
\draw  [fill={rgb, 255:red, 0; green, 0; blue, 0 }  ,fill opacity=1 ] (315,105) .. controls (315,103.34) and (316.34,102) .. (318,102) .. controls (319.66,102) and (321,103.34) .. (321,105) .. controls (321,106.66) and (319.66,108) .. (318,108) .. controls (316.34,108) and (315,106.66) .. (315,105) -- cycle ;
\draw  [draw opacity=0][fill={rgb, 255:red, 155; green, 155; blue, 155 }  ,fill opacity=0.25 ] (201.8,20.6) -- (318,20.6) -- (318,105) -- (201.8,105) -- cycle ;
\draw  [draw opacity=0][fill={rgb, 255:red, 155; green, 155; blue, 155 }  ,fill opacity=0.25 ] (318,105) -- (431.4,105) -- (431.4,190.2) -- (318,190.2) -- cycle ;
\draw [color={rgb, 255:red, 208; green, 2; blue, 27 }  ,draw opacity=1 ]   (318,20.2) -- (318,102) ;

\draw (450,159.6) node [anchor=north west][inner sep=0.75pt]    {$f_{1}$};
\draw (224.8,7.2) node [anchor=north west][inner sep=0.75pt]    {$f_{2}$};
\draw (304.4,82.4) node [anchor=north west][inner sep=0.75pt]    {$\mathbf{f}( x)$};
\draw (322.8,24.6) node [anchor=north west][inner sep=0.75pt]  [font=\small] [align=left] {{\scriptsize Region dominated by $\displaystyle \mathbf{f}( x)$}};
\draw (216.4,113.4) node [anchor=north west][inner sep=0.75pt]  [font=\small] [align=left] {{\scriptsize Region dominates $\displaystyle \mathbf{f}( x)$}};
\draw (334.4,56.2) node [anchor=north west][inner sep=0.75pt]    {$\{\mathbf{f}( x) \ \prec \dotsc \}$};
\draw (218.4,143) node [anchor=north west][inner sep=0.75pt]    {$\{\dotsc \prec \mathbf{f}( x)\}$};
\draw (222.8,54.6) node [anchor=north west][inner sep=0.75pt]  [font=\footnotesize] [align=left] {"Incomparable"};
\draw (336,138.6) node [anchor=north west][inner sep=0.75pt]  [font=\footnotesize] [align=left] {"Incomparable"};

\end{tikzpicture}}
    \caption{Illustration of the different domination regions in Multiobjective optimization}
    \label{fig:domination_relation}
\end{figure}
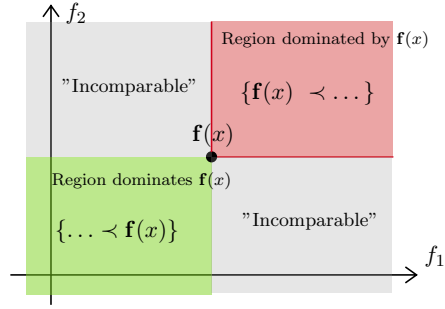

The Pareto front, which can be understood as the multiobjective counterpart to
the minimum, is defined as the set of all the non-dominated objectives
\begin{align}
\ParFront &= \{\fb(x)\text{ s.t. } x\in \Xspace\text{ and } \fb(x) \text{ non dominated } \} \\
&= \{\fb(x) \text{ s.t. } x \in \Xspace\text{ and } \nexists x'\in \Xspace, \fb(x') \prec \fb(x)\} \notag\,.
\end{align}

The Pareto set is the preimage of $\ParFront$, i.e. the set of Pareto-optimal points, noted $\ParSet$:
\begin{align}
\ParSet  = \{x\in\Xspace \text{ s.t. } \nexists x',  \fb(x') \prec  \fb(x)\}\,.
\end{align}
An example of a Pareto Front and Pareto Set is shown~\cref{fig:ill_pareto_front_set} with a discrete input space $\Xspace$.

\begin{figure}
    \centering
    \scalebox{0.8}{\tikzset{every picture/.style={line width=0.75pt}} 

\begin{tikzpicture}[x=0.75pt,y=0.75pt,yscale=-1,xscale=1]

\draw  (55,200.3) -- (231,200.3)(72.6,59) -- (72.6,216) (224,195.3) -- (231,200.3) -- (224,205.3) (67.6,66) -- (72.6,59) -- (77.6,66)  ;
\draw  (348,200.3) -- (551,200.3)(368.3,50) -- (368.3,217) (544,195.3) -- (551,200.3) -- (544,205.3) (363.3,57) -- (368.3,50) -- (373.3,57)  ;
\draw  [draw opacity=0][fill={rgb, 255:red, 74; green, 144; blue, 226 }  ,fill opacity=1 ] (84,160) .. controls (84,158.9) and (84.9,158) .. (86,158) .. controls (87.1,158) and (88,158.9) .. (88,160) .. controls (88,161.1) and (87.1,162) .. (86,162) .. controls (84.9,162) and (84,161.1) .. (84,160) -- cycle ;
\draw  [draw opacity=0][fill={rgb, 255:red, 208; green, 2; blue, 27 }  ,fill opacity=1 ] (113.5,75.5) .. controls (113.5,74.4) and (114.4,73.5) .. (115.5,73.5) .. controls (116.6,73.5) and (117.5,74.4) .. (117.5,75.5) .. controls (117.5,76.6) and (116.6,77.5) .. (115.5,77.5) .. controls (114.4,77.5) and (113.5,76.6) .. (113.5,75.5) -- cycle ;

\draw  [draw opacity=0][fill={rgb, 255:red, 126; green, 211; blue, 33 }  ,fill opacity=1 ] (182,158) .. controls (182,156.9) and (182.9,156) .. (184,156) .. controls (185.1,156) and (186,156.9) .. (186,158) .. controls (186,159.1) and (185.1,160) .. (184,160) .. controls (182.9,160) and (182,159.1) .. (182,158) -- cycle ;
\draw    (220,102.5) .. controls (259.6,72.8) and (312.92,70.54) .. (353.28,100.09) ;
\draw [shift={(354.5,101)}, rotate = 216.98] [color={rgb, 255:red, 0; green, 0; blue, 0 }  ][line width=0.75]    (10.93,-3.29) .. controls (6.95,-1.4) and (3.31,-0.3) .. (0,0) .. controls (3.31,0.3) and (6.95,1.4) .. (10.93,3.29)   ;
\draw  [draw opacity=0][fill={rgb, 255:red, 74; green, 144; blue, 226 }  ,fill opacity=1 ] (495.75,71.25) .. controls (495.75,70.15) and (496.65,69.25) .. (497.75,69.25) .. controls (498.85,69.25) and (499.75,70.15) .. (499.75,71.25) .. controls (499.75,72.35) and (498.85,73.25) .. (497.75,73.25) .. controls (496.65,73.25) and (495.75,72.35) .. (495.75,71.25) -- cycle ;
\draw [color={rgb, 255:red, 245; green, 166; blue, 35 }  ,draw opacity=1 ][line width=1.5]    (395,38) -- (395,101) ;
\draw [color={rgb, 255:red, 245; green, 166; blue, 35 }  ,draw opacity=1 ][line width=1.5]    (458,101) -- (458,165.5) ;
\draw [color={rgb, 255:red, 245; green, 166; blue, 35 }  ,draw opacity=1 ][line width=1.5]    (458,101) -- (395,101) ;
\draw [color={rgb, 255:red, 245; green, 166; blue, 35 }  ,draw opacity=1 ][line width=1.5]    (568.5,165.5) -- (458,165.5) ;

\draw  [draw opacity=0][fill={rgb, 255:red, 245; green, 166; blue, 35 }  ,fill opacity=0.24 ] (395.8,39.6) -- (571.05,39.6) -- (571.05,102.85) -- (395.8,102.85) -- cycle ;
\draw  [draw opacity=0][fill={rgb, 255:red, 245; green, 166; blue, 35 }  ,fill opacity=0.24 ] (458.8,102.6) -- (570.8,102.6) -- (570.8,165.85) -- (458.8,165.85) -- cycle ;

\draw  [draw opacity=0][fill={rgb, 255:red, 208; green, 2; blue, 27 }  ,fill opacity=1 ] (393,101) .. controls (393,99.9) and (393.9,99) .. (395,99) .. controls (396.1,99) and (397,99.9) .. (397,101) .. controls (397,102.1) and (396.1,103) .. (395,103) .. controls (393.9,103) and (393,102.1) .. (393,101) -- cycle ;
\draw  [draw opacity=0][fill={rgb, 255:red, 126; green, 211; blue, 33 }  ,fill opacity=1 ] (456,165.5) .. controls (456,164.4) and (456.9,163.5) .. (458,163.5) .. controls (459.1,163.5) and (460,164.4) .. (460,165.5) .. controls (460,166.6) and (459.1,167.5) .. (458,167.5) .. controls (456.9,167.5) and (456,166.6) .. (456,165.5) -- cycle ;
\draw  [draw opacity=0][fill={rgb, 255:red, 189; green, 16; blue, 224 }  ,fill opacity=1 ] (133.5,145.5) .. controls (133.5,144.4) and (134.4,143.5) .. (135.5,143.5) .. controls (136.6,143.5) and (137.5,144.4) .. (137.5,145.5) .. controls (137.5,146.6) and (136.6,147.5) .. (135.5,147.5) .. controls (134.4,147.5) and (133.5,146.6) .. (133.5,145.5) -- cycle ;
\draw  [draw opacity=0][fill={rgb, 255:red, 189; green, 16; blue, 224 }  ,fill opacity=1 ] (393,61.5) .. controls (393,60.4) and (393.9,59.5) .. (395,59.5) .. controls (396.1,59.5) and (397,60.4) .. (397,61.5) .. controls (397,62.6) and (396.1,63.5) .. (395,63.5) .. controls (393.9,63.5) and (393,62.6) .. (393,61.5) -- cycle ;

\draw (277.5,56.9) node [anchor=north west][inner sep=0.75pt]  [font=\large]  {$\mathbf{f}$};
\draw (92.5,11) node [anchor=north west][inner sep=0.75pt]   [align=left] {Control Space};
\draw (402.5,11.5) node [anchor=north west][inner sep=0.75pt]   [align=left] {Objective Space};
\draw (551.5,189.9) node [anchor=north west][inner sep=0.75pt]    {$f_{1}$};
\draw (358,29.4) node [anchor=north west][inner sep=0.75pt]    {$f_{2}$};
\draw (118,57.4) node [anchor=north west][inner sep=0.75pt]  [color={rgb, 255:red, 255; green, 0; blue, 0 }  ,opacity=1 ]  {$x_{1}$};
\draw (401.6,75.51) node [anchor=north west][inner sep=0.75pt]  [font=\footnotesize,color={rgb, 255:red, 208; green, 2; blue, 27 }  ,opacity=1 ]  {$\mathbf{f}( x_{1})$};
\draw (470,173.9) node [anchor=north west][inner sep=0.75pt]  [font=\footnotesize,color={rgb, 255:red, 65; green, 117; blue, 5 }  ,opacity=1 ]  {$\mathbf{f}( x_{2})$};
\draw (489.5,46.4) node [anchor=north west][inner sep=0.75pt]  [font=\footnotesize,color={rgb, 255:red, 74; green, 144; blue, 226 }  ,opacity=1 ]  {$\mathbf{f}( x_{3})$};
\draw (186.5,135.4) node [anchor=north west][inner sep=0.75pt]  [color={rgb, 255:red, 65; green, 117; blue, 5 }  ,opacity=1 ]  {$x_{2}$};
\draw (86.5,139.4) node [anchor=north west][inner sep=0.75pt]  [color={rgb, 255:red, 74; green, 144; blue, 226 }  ,opacity=1 ]  {$x_{3}$};
\draw (468.03,93.27) node [anchor=north west][inner sep=0.75pt]  [font=\scriptsize,color={rgb, 255:red, 245; green, 102; blue, 35 }  ,opacity=1 ] [align=left] {\begin{minipage}[lt]{60.63pt}\setlength\topsep{0pt}
Region dominated
\begin{center}
by $\displaystyle f(\mathcal{X})$
\end{center}

\end{minipage}};
\draw (565.6,164.8) node [anchor=north west][inner sep=0.75pt]  [color={rgb, 255:red, 245; green, 132; blue, 35 }  ,opacity=1 ]  {$\mathcal{P}^{*}$};
\draw (136,124.9) node [anchor=north west][inner sep=0.75pt]  [color={rgb, 255:red, 189; green, 16; blue, 224 }  ,opacity=1 ]  {$x_{4}$};
\draw (401,42.4) node [anchor=north west][inner sep=0.75pt]  [font=\footnotesize,color={rgb, 255:red, 74; green, 144; blue, 226 }  ,opacity=1 ]  {$\textcolor[rgb]{0.74,0.06,0.88}{\mathbf{f}( x_{4})}$};

\end{tikzpicture}}
    \caption{Pareto front for $\Xspace = \{x_1,x_2,x_3, x_4\}$. In this case, $\ParFront = \{\fb(x_1),\fb(x_2)\}$ (since $\fb(x_1) \prec \fb(x_4)$) and $\ParSet = \{x_1,x_2\}$.}
    \label{fig:ill_pareto_front_set}
\end{figure}
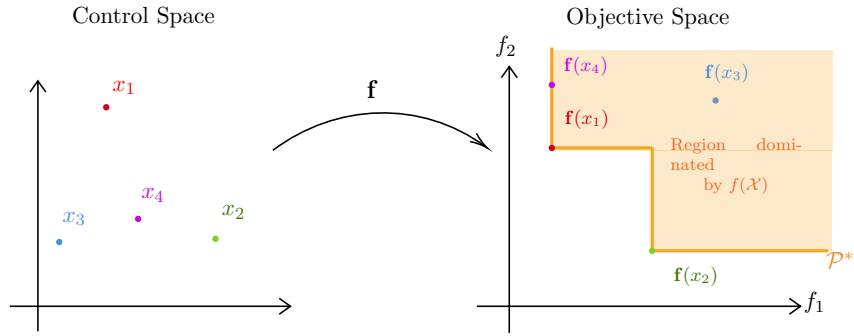

In practice, solving the multiobjective optimization problem means finding a finite approximation of the
Pareto Set and its associated Pareto Front, so that the decision maker can make
its choice among all the best possible \vicreplace{compromises}{tradeoffs}.

Pareto optimal solutions are typically approximated by successive scalarizations
of the multiobjective problem \cite{miettinen_nonlinear_1998,zhang2020random}
and application of a mono-objective optimization algorithm, or by adaptive
stochastic sampling algorithms such as the multiobjective version of CMA-ES
\cite{igel_covariance_2007,toure2019uncrowded} or NSGA-II \cite{deb_fast_2002}.
In the case of computationally expensive objectives, several methods based on
Bayesian Optimization have been derived that bring
together the building of surrogates to the true functions and the optimization.
We present them succinctly in what follows.

\subsection{Bayesian Multiobjective Optimization}
\label{ssec:bayesian_MO}
In a similar fashion as in single objective Bayesian Optimization, we can model
the objective function using Gaussian Processes (see
\cite{rasmussen_gaussian_2006,shahriari_taking_2016} for introductions).
We assume that the different objectives are independent, and use a zero-mean Gaussian process with Matérn 5/2 kernel as prior, denoted as $\Fb$}.

Let $\mathcal{D} = \left\{(x_1, \fb(x_1)),\dots,(x_{n_0}, \fb(x_{n_0})\right\}$
be the initial design of experiment, that is the set of already evaluated
input-output pairs. By conditioning $\Fb$ on the design $\mathcal{D}$, we obtain
\begin{equation}
    \Fb \mid \mathcal{D}\sim \mathrm{GP}(\mb_\Fb, \mathbf{k}_{\Fb})\,,
    \label{eq:gp_x}
\end{equation}
where $\mb_\Fb: \Xspace \rightarrow \mathbb{R}^d$ is called the GP prediction,
GP mean or kriging mean, and $\mathbf{k}_\Fb: \Xspace \times \Xspace \rightarrow
\mathbb{R}^{d \times d}$ is the covariance function. \vicdelete{, which is based on a
parametric kernel (Gaussian, Matérn 3/5 or Matérn 5/2 in most cases) for classical GP
regression.}\vicadd{Since the objective are independent, $\mathbf{k}_\Fb(x,x)$ is a diagonal matrix for any $x\in\Xspace$.}

The conditioning on the design of experiment will be dropped if it
is clear from the context. At a specific input $x\in\Xspace $, we have by
properties of the GP

\begin{equation}
    \Fb(x) \sim \mathcal{N}(\mb_\Fb(x), \mathbf{k}_{\Fb}(x,x))\,, \quad
\end{equation}
\vicdelete{where $\mathbf{S}^2_{\Fb}(x) = \mathbf{k}_{\Fb}\left(x, x\right) \in
\mathbb{R}^{d\times d}$. }This surrogate, which takes into account both a
prediction $\mb_{\Fb}$, and a measure of the uncertainty associated with this
prediction through the kernel $\mathbf{k}_{\Fb}$.
Based on this probabilistic representation of the unknown function, we can
define a progress measure $\alpha:\Xspace \rightarrow \mathbb{R}$ on the input
space, also known as acquisition function
\cite{frazier2018tutorial,garnett_bayesian_2023}. It quantifies the interest we
have at evaluating a specific choice of the control variables by the true
function. This forms the core of the typical Bayesian Optimization loop, as
described in the pseudocode of \cref{alg:bo_loop}. 
\begin{algorithm}
\begin{algorithmic}
\Require GP prior $\Fb$, Initial Design of Experiment $\mathcal{D}$
\While{Budget not exceeded}
\State Condition $\Fb$ on the current DoE $\mathcal{D}$ 
\State Optimize the acquisition function $\alpha$ to get $x_{\text{next}}$ 
\State Evaluate the true function $\fb(x_{\text{next}})$
\State Update DoE: $\mathcal{D}\gets\mathcal{D} \cup (x_{\text{next}}, \fb(x_{\text{next}}))$ 
\EndWhile
\end{algorithmic}
\caption{Bayesian Optimization Loop}
\label{alg:bo_loop}
\end{algorithm}

Bayesian Optimization has also been applied to Multiobjective Optimization,
where many progress measures have been defined and successfully implemented,
which rely on different properties and characterizations of the Pareto Front and
Set. \vicadd{A thorough survey can be found in \cite{rojas-gonzalez_survey_2020}.} Some criteria, such as the ParEGO found in
\cite{knowles_parego_2006} rely on the \emph{scalarization} of the objective
vector, i.e., on the aggregation of all the objectives into a single scalar
(through convex combination for instance) that is optimized afterwards. Pareto
Active Learning (PAL), as described in~\cite{zuluaga_active_2013} is based on
the classification of points as Pareto optimal or
not. 
In \cite{picheny_multiobjective_2015}, the author proposes a
measure of ``global'' uncertainty, which measures the uncertainty on the Pareto
front as the integral (in the control space) of the probability of improvement.
This measure of uncertainty can then be optimized using a SUR method
(\cite{bect_supermartingale_2019}). Another type of uncertainty measures rely
on an information-theoretic approach as done in
\cite{hernandez-lobato_predictive_2016,belakaria_maxvalue_2020,tu_joint_2022}.

One of the most used criterion is the Expected Hypervolume Improvement (EHVI)
which can be thought of as a natural extension of the well-known EI criterion of
\cite{schonlau1998global,jones_efficient_1998}. This popular criterion
introduced in~\cite{emmerich_single_2006} still presents some computational
challenges such as efficient partitioning of the dominated space which have been
further studied
in~\cite{ponweiser_multiobjective_2008,yang_multiobjective_2019,daulton_differentiable_2020,
daulton_parallel_2021}.
This criterion relies on the hypervolume of the dominated region, which
possesses monotonicity properties with respect to the domination
relation~\cite{audet_performance_2021}. In order to upper-bound this volume, it
is necessary to introduce $\mathbb{B}_{\text{ref}} = \{ \yb \mid \yb \prec
\yb_{\text{ref}}\}$ where $\yb_{\text{ref}}$ is a reference point which is chosen usually as
dominated by the Nadir point. Let $\hat{\ParFront}$ be an approximation of the
Pareto front obtained for instance using the GP prediction or using already
evaluated points. It represents the current best approximation of the Pareto front. The hypervolume of the region dominated by $\hat{\ParFront}$
is defined as
\begin{equation}
    \HV(\hat{\ParFront}) = \int_{\mathbb{B}_{\text{ref}}} \mathbbm{1}_{\{\hat{\ParFront} \prec \yb\}} \, \mathrm{d}\yb \,.
\end{equation}
\vicadd{where $\mathbbm{1}_{\{\hat{\ParFront} \prec \yb\}} = 1$ if $\yb$ is dominated by $\ParFront$, $0$ elsewhere.}
The improvement in the hypervolume of the region dominated by $\hat{\ParFront}$
when adding $\yb$ to the approximation can be written as
\begin{equation}
    \mathrm{HVI}\left(\yb, \hat{\ParFront}\right) = \HV(\hat{\ParFront}\cup \{\yb\}) - \HV(\hat{\ParFront}) \,. 
    \label{eq:HVI}
\end{equation}

Since $\Fb$ is modeled using a GP, we have that the future evaluation $\yb$ is
distributed according to $\Fb(x)$, which is multivariate normal. Averaging this
improvement with respect to $\Fb(x)$ yields the Expected Hypervolume
Improvement:
\begin{equation}
    \mathop{\mathrm{EHVI}}(x) = \Ex_{\Fb(x)}\left[\mathrm{HVI}(\Fb(x), \hat{\ParFront}) \right] \,.
    \label{eq:EHVI}
\end{equation}
Variations \vicadd{around}\vicdelete{of} this acquisition function have been studied, where the reference
point for the computation of the dominated region is chosen specifically to
target some regions of the Pareto front, as in \cite{gaudrie_targeting_2020}.
Instead of using the hypervolume, the maximin marginal improvement can be
considered as in \cite{balling_maximin_2003,bautista_sequential_2009,svenson_multiobjective_2016}.

\section{Multiobjective Optimization in the presence of uncertainties and related work}
\label{sec:moouu}
We consider now that the different objectives  are functions of the control
variable $x$ and another variable $u\in \Uspace\subseteq \mathbb{R}^{n_U}$:
\begin{equation}
\begin{array}{rcl}
\fb: \Xspace \times \Uspace & \longrightarrow& \mathbb{R}^d \\
(x, u) & \longmapsto &\fb(x, u) = (f_1(x, u),\dots,f_d(x, u))
\end{array}
\end{equation}
This additional variable $u$ represents some uncertainties in the optimization
problem due to environmental conditions, which are either uncontrollable, or unknown to the
modeler when the choice of $x$ must be made. This formalism helps tackle the problem
in a large number of situations as in~\cite{rivier_surrogateassisted_2022} and
\cite{inatsu_bounding_2023}.
We assume in the following that this environmental variable is modeled using a
random variable $U$ of known probability density function $p_U$, with support
included in $\Uspace$. This modeling assumes then that the function $\fb$ is
deterministic, so that sampling $U$ is up to the user. This is sometimes
named a \emph{simulator} setting for $u$, as opposed to \emph{uncontrollable}
setting, where $u$ cannot be controlled even during the optimization (with a
stochastic simulator for instance).

Multiobjective optimization under uncertainties as been treated in various ways
in the literature. In~\cite{tu_scalarisationbased_2025}, the authors review
scalarization-based methods for robust multiobjective optimization, while in
\cite{daulton_robust_2022}, the authors propose to use the Multivariate
Value-at-Risk (MVaR), which can be seen as a multidimensional extension of the
Value-at-Risk, which is then optimized using a scalarization method. Working directly on some statistics of the objective function is also popular:
in \cite{rivier_surrogateassisted_2022} and \cite{inatsu_bounding_2023}, the
authors consider some statistics of $\fb(x, U)$ which are then optimized. \vicadd{This approach has been also well studied when considering stochastic simulators, as in \cite{rojasgonzalez_multiobjective_2020,daulton_parallel_2021,pal_multiobjective_2020}}.
The case of the uncertainties appearing as random perturbations of the control
variable has been treated more specifically in
\cite{gutjahr_stochastic_2016,peitz_survey_2018,ribaud_robust_2020}.

One natural approach is to look for the Pareto front and Pareto set of the
expected value of the objective vector, that is
\begin{equation}
    \min_{x \in \Xspace} \Ex_U[\fb(x, U)] \,.
    \label{eq:mean_obj_pareto}
\end{equation}

\begin{figure}[ht]
    \centering
    \scalebox{0.8}{\tikzset{every picture/.style={line width=0.75pt}} 

\begin{tikzpicture}[x=0.75pt,y=0.75pt,yscale=-1,xscale=1]

\draw  (96.33,170.57) -- (232.33,170.57)(109.93,61.67) -- (109.93,182.67) (225.33,165.57) -- (232.33,170.57) -- (225.33,175.57) (104.93,68.67) -- (109.93,61.67) -- (114.93,68.67)  ;
\draw  (257,169.9) -- (393,169.9)(270.6,61) -- (270.6,182) (386,164.9) -- (393,169.9) -- (386,174.9) (265.6,68) -- (270.6,61) -- (275.6,68)  ;
\draw  [color={rgb, 255:red, 208; green, 2; blue, 27 }  ,draw opacity=1 ] (135.52,93.03) .. controls (138.62,88.91) and (156.73,97.32) .. (175.97,111.81) .. controls (195.2,126.3) and (208.28,141.38) .. (205.18,145.5) .. controls (202.07,149.63) and (183.96,141.22) .. (164.73,126.73) .. controls (145.49,112.24) and (132.42,97.16) .. (135.52,93.03) -- cycle ;
\draw  [color={rgb, 255:red, 208; green, 2; blue, 27 }  ,draw opacity=1 ] (148.94,103.14) .. controls (150.85,100.61) and (161.98,105.77) .. (173.8,114.68) .. controls (185.63,123.59) and (193.67,132.86) .. (191.76,135.4) .. controls (189.85,137.93) and (178.72,132.76) .. (166.89,123.86) .. controls (155.07,114.95) and (147.03,105.67) .. (148.94,103.14) -- cycle ;
\draw  [color={rgb, 255:red, 208; green, 2; blue, 27 }  ,draw opacity=1 ] (128.96,88.09) .. controls (133.78,81.69) and (156.21,90.47) .. (179.07,107.69) .. controls (201.93,124.9) and (216.56,144.05) .. (211.74,150.45) .. controls (206.92,156.84) and (184.48,148.07) .. (161.62,130.85) .. controls (138.76,113.63) and (124.14,94.49) .. (128.96,88.09) -- cycle ;

\draw  [color={rgb, 255:red, 208; green, 2; blue, 27 }  ,draw opacity=1 ] (364.24,92.77) .. controls (367.54,96.74) and (355.18,112.42) .. (336.65,127.79) .. controls (318.11,143.17) and (300.41,152.41) .. (297.12,148.44) .. controls (293.82,144.47) and (306.18,128.78) .. (324.72,113.41) .. controls (343.25,98.04) and (360.95,88.79) .. (364.24,92.77) -- cycle ;
\draw  [color={rgb, 255:red, 208; green, 2; blue, 27 }  ,draw opacity=1 ] (351.31,103.49) .. controls (353.34,105.93) and (345.74,115.57) .. (334.35,125.02) .. controls (322.95,134.47) and (312.07,140.16) .. (310.05,137.71) .. controls (308.02,135.27) and (315.62,125.63) .. (327.01,116.18) .. controls (338.41,106.73) and (349.29,101.05) .. (351.31,103.49) -- cycle ;
\draw  [color={rgb, 255:red, 208; green, 2; blue, 27 }  ,draw opacity=1 ] (370.57,87.52) .. controls (375.68,93.69) and (361.97,113.5) .. (339.94,131.76) .. controls (317.91,150.03) and (295.91,159.85) .. (290.8,153.68) .. controls (285.68,147.52) and (299.4,127.71) .. (321.42,109.44) .. controls (343.45,91.17) and (365.45,81.36) .. (370.57,87.52) -- cycle ;
\draw  [draw opacity=0][fill={rgb, 255:red, 208; green, 2; blue, 27 }  ,fill opacity=1 ] (168.18,119.27) .. controls (168.18,118.07) and (169.15,117.1) .. (170.35,117.1) .. controls (171.54,117.1) and (172.51,118.07) .. (172.51,119.27) .. controls (172.51,120.47) and (171.54,121.44) .. (170.35,121.44) .. controls (169.15,121.44) and (168.18,120.47) .. (168.18,119.27) -- cycle ;
\draw  [draw opacity=0][fill={rgb, 255:red, 208; green, 2; blue, 27 }  ,fill opacity=1 ] (328.51,120.6) .. controls (328.51,119.41) and (329.48,118.44) .. (330.68,118.44) .. controls (331.88,118.44) and (332.85,119.41) .. (332.85,120.6) .. controls (332.85,121.8) and (331.88,122.77) .. (330.68,122.77) .. controls (329.48,122.77) and (328.51,121.8) .. (328.51,120.6) -- cycle ;
\draw  [draw opacity=0][fill={rgb, 255:red, 208; green, 2; blue, 27 }  ,fill opacity=0.33 ] (170.33,184.67) .. controls (156.67,184.67) and (169.67,170.67) .. (109.93,170.57) .. controls (50.2,170.47) and (285.67,171.33) .. (229.67,170.67) .. controls (173.67,170) and (184,184.67) .. (170.33,184.67) -- cycle ;
\draw  [draw opacity=0][fill={rgb, 255:red, 208; green, 2; blue, 27 }  ,fill opacity=0.33 ] (95.86,121.36) .. controls (95.86,107.7) and (109.86,120.7) .. (109.96,60.96) .. controls (110.06,1.23) and (109.19,236.7) .. (109.86,180.7) .. controls (110.53,124.7) and (95.86,135.03) .. (95.86,121.36) -- cycle ;
\draw  [draw opacity=0][fill={rgb, 255:red, 208; green, 2; blue, 27 }  ,fill opacity=0.33 ] (331,184) .. controls (317.33,184) and (330.33,170) .. (270.6,169.9) .. controls (210.87,169.8) and (446.33,170.67) .. (390.33,170) .. controls (334.33,169.33) and (344.67,184) .. (331,184) -- cycle ;
\draw  [draw opacity=0][fill={rgb, 255:red, 208; green, 2; blue, 27 }  ,fill opacity=0.33 ] (256.53,122.03) .. controls (256.53,108.36) and (270.53,121.36) .. (270.63,61.63) .. controls (270.73,1.9) and (269.86,237.36) .. (270.53,181.36) .. controls (271.19,125.36) and (256.53,135.7) .. (256.53,122.03) -- cycle ;

\draw (148,65.73) node [anchor=north west][inner sep=0.75pt]    {$\mathbf{f}( x_{1} ,U)$};
\draw (303.33,63.07) node [anchor=north west][inner sep=0.75pt]    {$\mathbf{f}( x_{2} ,U)$};
\draw (226.67,176.07) node [anchor=north west][inner sep=0.75pt]    {$f_{1}$};
\draw (392,169.4) node [anchor=north west][inner sep=0.75pt]    {$f_{1}$};
\draw (279,47.4) node [anchor=north west][inner sep=0.75pt]    {$f_{2}$};
\draw (87.67,48.07) node [anchor=north west][inner sep=0.75pt]    {$f_{2}$};

\end{tikzpicture}}
    \caption{Sketch of the different behaviors: Both $x_1$ and $x_2$ provide the same mean objective, the marginals of the objective vector are the same, but the behavior in a multiobjective optimization problem is not the same}
    \label{fig:marginals}
\end{figure}
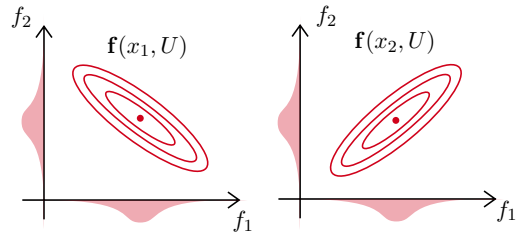

Such a formulation allows to remove the uncertainty from each component independently.
Similarly to the single objective case, mean minimization does not take into
account the variability of the solution, nor its skewness. In a multiobjective
setting, it fails also to consider the correlation between the objectives.
Indeed, if we take for instance a problem with two objectives, a positive
correlation between those indicates that under the uncertainty, both objectives
have the tendency to be degraded or improved simultaneously. If the correlation
is negative, the samples provided are more likely to be non-dominated
as progress in one of the objective often occurs while the other
objective is degraded. This is illustrated in \cref{fig:marginals}. 

Other risk measures have been defined for random vectors such as probabilistic
definitions of dominance \cite{rivier_surrogateassisted_2022,khosravi_probabilistic_2018} or
\cite{ek_learning_2021}. \vicadd{These works usually introduce an additional parameter, which acts as a threshold to assess whether a design dominates another one with a large enough probability.}

\vicreplace{In this work, we will instead directly consider the solutions of the
multiobjective optimization problem, namely the Pareto front and Pareto set, as
random quantities.}{The main novelty of this work is to consider the solutions of the multiobjective problem, namely the Pareto front and Pareto set, as random quantities, and to exploit some statistical information on those. In particular, we are interested in the probability of coverage, which indicates the frequency of which we can expect to be close to the optimal solutions. This allows to rank easily solutions in an interpretable way. We also propose a new method based on Bayesian Optimization in order to improve the estimation of the Conditional Pareto fronts and Conditional Pareto Sets.}

\section{Conditional Pareto front, robustness through the probability of coverage}
\label{sec:bmoouu}
\subsection{Definition and probability of coverage}
Let us consider a sample $u\in\Uspace$ fixed. The conditional multiobjective
optimization problem associated with $x\mapsto\fb(x, u)$ is
\begin{align}
    \min_{x \in \Xspace} \fb(x, u) \,\, ,
\end{align}
which is fully deterministic.
The solution of this MOO problem is the Conditional Pareto Front (CPF)
$\ParFront(u)$, and the Conditional Pareto Set (CPS) $\ParSet(u)$.
The idea of looking at the CPF and CPS can also be
found in \cite{ide_relationship_2014,ide_robustness_2016} in a non-probabilistic setting.

\begin{figure}[ht]
    \centering
    \scalebox{0.9}{\input{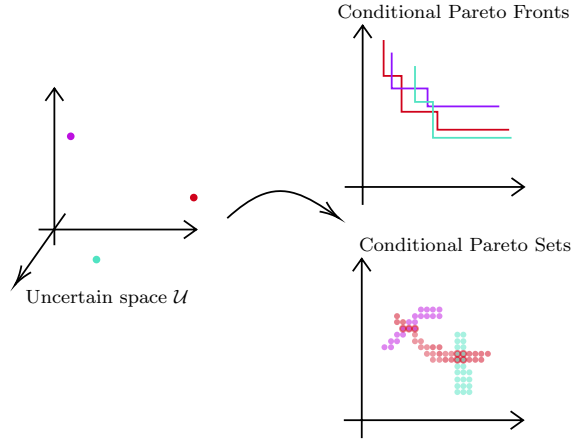}}
    \caption{Illustration of {discretized} Conditional Pareto Sets and Fronts, depending on the environmental variable }
    \label{fig:schema_cpf}
\end{figure}

We are interested in the coverage probability of the random closed set
$\ParSet(U)$ that is, the probability that a given $x\in\Xspace$ belongs to the CPS,
\begin{align}
    \Prob_U\left[x \in \ParSet(U)\right] = \Ex_U\left[\mathbbm{1}_{\ParSet(U)}(x)\right]\,.
    \label{eq:prob_coverage}
\end{align}

\begin{figure}
    \centering
    \includegraphics[width=0.9\linewidth]{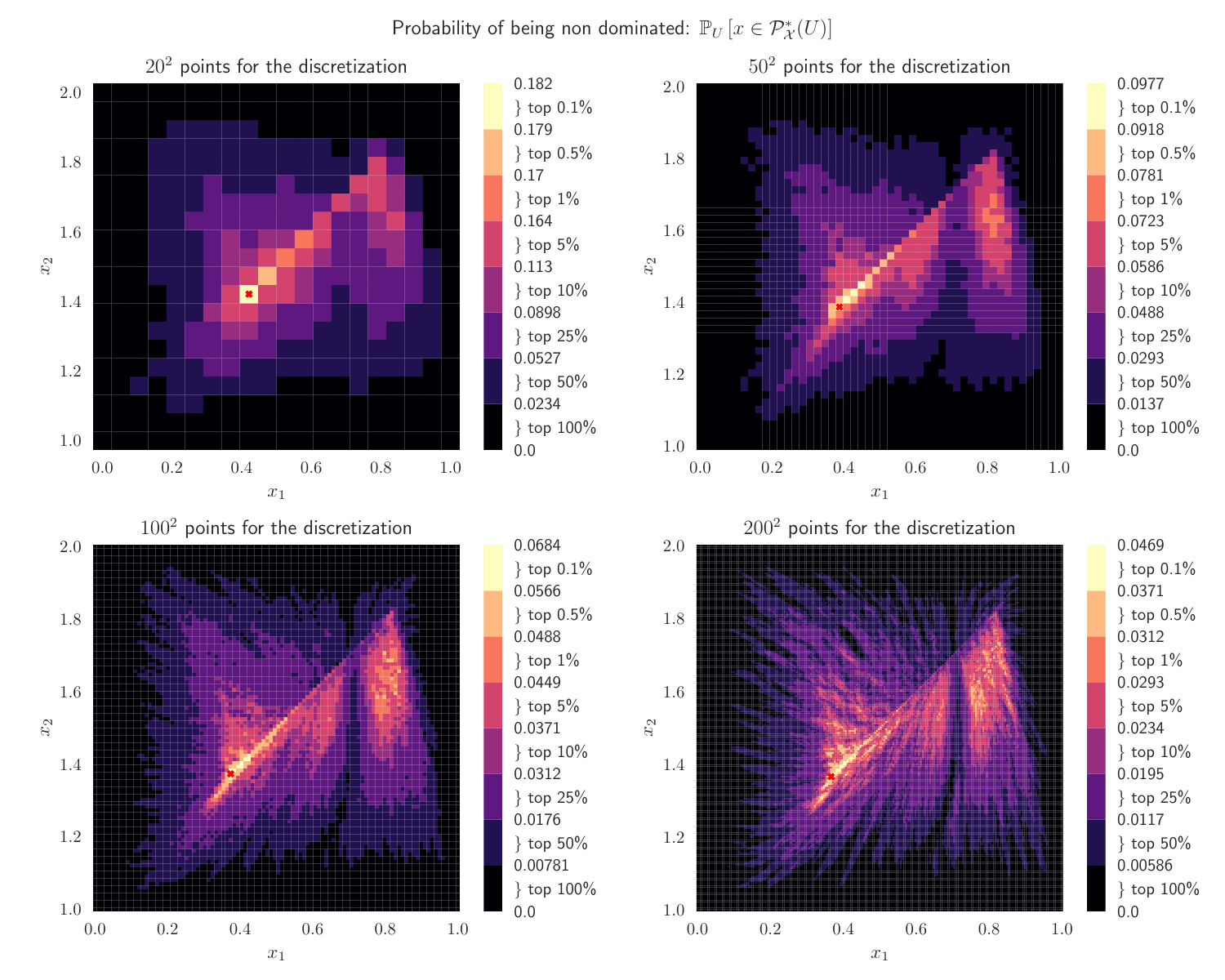}
    \caption{Probability of coverage of the Pareto Set for the problem described
    (\ref{eq:def_f2x2}) \vicadd{with different discretizations. The levels are chosen to represent specific percentiles of values}. \vicdelete{Best 5\% of values indicated using the green line,
    best 1\% by the red one.} Maximum indicated with the red cross \vicadd{ at approximately $(0.4, 1.4)$}.}
    \label{fig:prob_cpf}
\end{figure}

The CPS \vicreplace{can be however of measure $0$}{has however zero measure in most cases}. We rely then on a \vicreplace{discretization}{discretized} version
of the problem, where the CPS is computed among a large set of candidate points.
\vicadd{\subsection{Approximation based on the discretization}
Let us consider $\gb = \fb(\cdot,u)$ for a $u\in\Uspace$.
We define the $\epsilon$-accurate Pareto front as 
\begin{equation}
    \ParFront_{\epsilon} =  \{\gb(x) \text{ s.t. } x \in \Xspace\text{ and } \nexists x'\in \Xspace,\quad \gb(x') \prec \gb(x)-\epsilon\} \notag\,.
\end{equation}
where $\gb(x)-\epsilon = (g_1(x)-\epsilon,\dots,g_2(x) - \epsilon)$ represent each objective "improved" by $\epsilon$.
We assume that we have a set $\mathfrak{X} \subset \Xspace^p$ of points, and the associated evaluations $\gb(\mathfrak{X)} = \{\gb(x) \text{ for } x\in \mathfrak{X}\}$. We denote the maximin distance as $\delta$, which is the maximum distance between a point of the domain, and a point of the design:
\begin{equation}
    \max_{x\in\Xspace}\min_{x_i \in \mathfrak{X}}\|x_i - x\|= \delta\,.
\end{equation}
Finally, we assume that the function is $L$-Lipschitz: $\|\gb(x_1) - \gb(x_2) \| \leq L \|x_1 - x_2 \|$.

Let $x^* \in \ParSet$. Since $\ParSet\subset\Xspace$, there exists $i$ such that $\|x_i - x^*\| \leq \delta$, meaning that any point in the Pareto set is at distance at most $\delta$ from an evaluated point.
The difference between the function value can then be bounded:
\begin{equation}
    \|\gb(x_i) - \gb(x^*) \| \leq L\|x_i - x^*\| \leq L \delta
\end{equation}
In other words, every point of the Pareto front is at most at a distance $L\delta$ from a point from the discretization, thus the set of non dominated points belongs to the $L\delta$-accurate Pareto front approximation. An illustration is provided on a numerical example \cref{fig:prob_cpf}. For simplicity, the probability of coverage in what follows is to be understood as its discretized version: a point belongs to the CPS when its evaluation is at a controlled distance to the CPF, distance which depends on the Lipschitz constant of the function, and the maximin distance of the discretization.  
}

\vicadd{In a context of multiobjective optimization, we can then}\vicdelete{We can} look for the point which maximizes the coverage
probability (\ref{eq:prob_coverage}) \vicadd{ or at least points with a large enough probability.} \vicdelete{Such a maximizer can also be interpreted as the Bayesian optimal
decision under the 0--1 loss.}
\vicdelete{ shows that for a $u \in \Uspace$, the problem
is in fact one of classification, where we need to evaluate whether $x$ belongs
to the conditional Pareto set $\ParSet(u)$.}
However, this coverage probability may be very expensive to compute with
sufficient precision. Indeed, a naive implementation would involve the
determination of Pareto sets for a sufficiently large number of samples of $U$.
If a finite number of candidates is considered 
multiobjective optimization procedure is already possibly expensive in itself,
so repeating this procedure for many different $u$ will quickly become
impossible. We propose then to use a surrogate model in order to alleviate the
number of evaluations of $\fb$, \vicreplace{an idea which belongs to the grand scheme of
Bayesian Optimization.}{an idea implemented within the context of Bayesian Optimization.}

\section{Bayesian Optimization for Multiobjective Optimization under Uncertainties}
Similarly as in (\ref{eq:gp_x}), we model the objective function using a GP
on the joint space $\Xspace\times \Uspace$. Let $\mathcal{D} =
\left\{((x_1,u_1), \fb(x_1,u_1)),\dots,((x_{n_0},u_{n_0}),
\fb(x_{n_0},u_{n_0}))\right\}$ be the initial design of experiment, and keeping
the same notation, we have
\begin{equation}
    \Fb \mid \mathcal{D}\sim \mathrm{GP}(\mb_\Fb, \mathbf{k}_{\Fb})\,,
\end{equation}
where $\mb_\Fb: \Xspace \times\Uspace \rightarrow \mathbb{R}^d$ is the kriging
mean, and $\mathbf{k}_\Fb: (\Xspace \times \Uspace) \times (\Xspace \times
\Uspace) \rightarrow \mathbb{R}^{d \times d}$ is the covariance function. \vicreplace{We can
use the kriging mean as a surrogate to estimate the probability of coverage.}{Based on this modelling of the unknown function, we can use the kriging mean as a plug-in to estimate the probability of coverage. This estimation can also be more or less conservative using the information coming from the covariance function.}


\subsection{Profile EHVI for CPF}
\label{sec:pehvi}
\begin{figure}[ht]
    \centering
    \scalebox{.8}{\tikzset{every picture/.style={line width=0.75pt}} 

\begin{tikzpicture}[x=0.75pt,y=0.75pt,yscale=-1,xscale=1]

\draw  [draw opacity=0][fill={rgb, 255:red, 80; green, 227; blue, 194 }  ,fill opacity=0.75 ] (498.89,115.06) -- (549.38,114.92) -- (549.38,127.63) -- (579.8,129.6) -- (579.8,140) -- (498.6,140) -- cycle ;
\draw [color={rgb, 255:red, 208; green, 2; blue, 27 }  ,draw opacity=1 ][line width=1.5]    (434,50.4) -- (434,115.2) -- (481.26,114.92) -- (538.6,114.92) -- (549.38,114.92) -- (549.38,129.6) -- (579.8,129.6) ;
\draw  (359,190.12) -- (586.77,190.12)(381.78,45.37) -- (381.78,206.2) (579.77,185.12) -- (586.77,190.12) -- (579.77,195.12) (376.78,52.37) -- (381.78,45.37) -- (386.78,52.37)  ;
\draw   (496.63,140) .. controls (496.63,138.91) and (497.51,138.03) .. (498.6,138.03) .. controls (499.69,138.03) and (500.57,138.91) .. (500.57,140) .. controls (500.57,141.09) and (499.69,141.97) .. (498.6,141.97) .. controls (497.51,141.97) and (496.63,141.09) .. (496.63,140) -- cycle ;
\draw  [fill={rgb, 255:red, 65; green, 117; blue, 5 }  ,fill opacity=1 ] (547.41,129.6) .. controls (547.41,128.51) and (548.3,127.63) .. (549.38,127.63) .. controls (550.47,127.63) and (551.35,128.51) .. (551.35,129.6) .. controls (551.35,130.69) and (550.47,131.57) .. (549.38,131.57) .. controls (548.3,131.57) and (547.41,130.69) .. (547.41,129.6) -- cycle ;
\draw   (578.63,50.4) .. controls (578.63,49.31) and (579.51,48.43) .. (580.6,48.43) .. controls (581.69,48.43) and (582.57,49.31) .. (582.57,50.4) .. controls (582.57,51.49) and (581.69,52.37) .. (580.6,52.37) .. controls (579.51,52.37) and (578.63,51.49) .. (578.63,50.4) -- cycle ;
\draw  [fill={rgb, 255:red, 65; green, 117; blue, 5 }  ,fill opacity=1 ] (432.03,115.2) .. controls (432.03,114.11) and (432.91,113.23) .. (434,113.23) .. controls (435.09,113.23) and (435.97,114.11) .. (435.97,115.2) .. controls (435.97,116.29) and (435.09,117.17) .. (434,117.17) .. controls (432.91,117.17) and (432.03,116.29) .. (432.03,115.2) -- cycle ;
\draw  [draw opacity=0][fill={rgb, 255:red, 189; green, 16; blue, 224 }  ,fill opacity=0.5 ] (451,115.2) -- (549.38,114.92) -- (549.38,120.4) -- (451,121.2) -- cycle ;
\draw   (449.03,121.2) .. controls (449.03,120.11) and (449.91,119.23) .. (451,119.23) .. controls (452.09,119.23) and (452.97,120.11) .. (452.97,121.2) .. controls (452.97,122.29) and (452.09,123.17) .. (451,123.17) .. controls (449.91,123.17) and (449.03,122.29) .. (449.03,121.2) -- cycle ;
\draw [color={rgb, 255:red, 208; green, 2; blue, 27 }  ,draw opacity=1 ][line width=1.5]    (109.12,47.42) -- (109.12,95.46) -- (156.37,95.46) -- (156.37,142.72) -- (224.5,142.72) -- (224.5,174.62) -- (253.01,174.62) ;
\draw  (58,190.12) -- (285.77,190.12)(80.78,45.37) -- (80.78,206.2) (278.77,185.12) -- (285.77,190.12) -- (278.77,195.12) (75.78,52.37) -- (80.78,45.37) -- (85.78,52.37)  ;
\draw  [draw opacity=0][fill={rgb, 255:red, 80; green, 227; blue, 194 }  ,fill opacity=0.75 ] (125.26,95.46) -- (156.37,95.46) -- (156.37,115.94) -- (125.26,115.94) -- cycle ;
\draw   (123.29,115.94) .. controls (123.29,114.85) and (124.17,113.97) .. (125.26,113.97) .. controls (126.35,113.97) and (127.23,114.85) .. (127.23,115.94) .. controls (127.23,117.03) and (126.35,117.91) .. (125.26,117.91) .. controls (124.17,117.91) and (123.29,117.03) .. (123.29,115.94) -- cycle ;
\draw  [draw opacity=0][fill={rgb, 255:red, 189; green, 16; blue, 224 }  ,fill opacity=0.5 ] (146.53,142.72) -- (224.5,142.72) -- (224.5,167.14) -- (146.53,167.14) -- cycle ;
\draw  [draw opacity=0][fill={rgb, 255:red, 189; green, 16; blue, 224 }  ,fill opacity=0.55 ] (146.53,95.46) -- (156.37,95.46) -- (156.37,142.72) -- (146.53,142.72) -- cycle ;
\draw   (144.29,167.14) .. controls (144.29,166.05) and (145.18,165.17) .. (146.26,165.17) .. controls (147.35,165.17) and (148.23,166.05) .. (148.23,167.14) .. controls (148.23,168.22) and (147.35,169.1) .. (146.26,169.1) .. controls (145.18,169.1) and (144.29,168.22) .. (144.29,167.14) -- cycle ;
\draw  [fill={rgb, 255:red, 65; green, 117; blue, 5 }  ,fill opacity=1 ] (107.15,95.46) .. controls (107.15,94.38) and (108.03,93.49) .. (109.12,93.49) .. controls (110.2,93.49) and (111.08,94.38) .. (111.08,95.46) .. controls (111.08,96.55) and (110.2,97.43) .. (109.12,97.43) .. controls (108.03,97.43) and (107.15,96.55) .. (107.15,95.46) -- cycle ;
\draw  [fill={rgb, 255:red, 65; green, 117; blue, 5 }  ,fill opacity=1 ] (154.4,142.72) .. controls (154.4,141.63) and (155.28,140.75) .. (156.37,140.75) .. controls (157.46,140.75) and (158.34,141.63) .. (158.34,142.72) .. controls (158.34,143.81) and (157.46,144.69) .. (156.37,144.69) .. controls (155.28,144.69) and (154.4,143.81) .. (154.4,142.72) -- cycle ;
\draw  [fill={rgb, 255:red, 65; green, 117; blue, 5 }  ,fill opacity=1 ] (222.53,174.62) .. controls (222.53,173.53) and (223.41,172.65) .. (224.5,172.65) .. controls (225.59,172.65) and (226.47,173.53) .. (226.47,174.62) .. controls (226.47,175.7) and (225.59,176.59) .. (224.5,176.59) .. controls (223.41,176.59) and (222.53,175.7) .. (222.53,174.62) -- cycle ;
\draw  [draw opacity=0][fill={rgb, 255:red, 155; green, 155; blue, 155 }  ,fill opacity=0.37 ] (109.04,51.36) -- (248.13,51.36) -- (248.13,94.99) -- (109.04,94.99) -- cycle ;
\draw  [draw opacity=0][fill={rgb, 255:red, 155; green, 155; blue, 155 }  ,fill opacity=0.37 ] (156.37,94.99) -- (248.13,94.99) -- (248.13,142.25) -- (156.37,142.25) -- cycle ;
\draw  [draw opacity=0][fill={rgb, 255:red, 155; green, 155; blue, 155 }  ,fill opacity=0.37 ] (224.97,142.25) -- (248.13,142.25) -- (248.13,174.62) -- (224.97,174.62) -- cycle ;
\draw   (246.16,51.36) .. controls (246.16,50.27) and (247.04,49.39) .. (248.13,49.39) .. controls (249.21,49.39) and (250.1,50.27) .. (250.1,51.36) .. controls (250.1,52.45) and (249.21,53.33) .. (248.13,53.33) .. controls (247.04,53.33) and (246.16,52.45) .. (246.16,51.36) -- cycle ;
\draw  [draw opacity=0][fill={rgb, 255:red, 155; green, 155; blue, 155 }  ,fill opacity=0.42 ] (580.6,50.4) -- (579.8,129.6) -- (549.38,129.6) -- (549.38,114.92) -- (434,115.2) -- (434,50.4) -- cycle ;

\draw (592.79,182.71) node [anchor=north west][inner sep=0.75pt]    {$f_{1}$};
\draw (359.74,37.95) node [anchor=north west][inner sep=0.75pt]    {$f_{2}$};
\draw (410.38,125.37) node [anchor=north west][inner sep=0.75pt]    {$\mathbf{F}( x_{2} ,u_{2})$};
\draw (582.92,31.95) node [anchor=north west][inner sep=0.75pt]   [align=left] {nadir};
\draw (406.4,23.2) node [anchor=north west][inner sep=0.75pt]  [color={rgb, 255:red, 208; green, 2; blue, 27 }  ,opacity=1 ]  {$\mathcal{P}^{*}_{\mathbf{m}}( u_{2})$};
\draw (484.8,144.75) node [anchor=north west][inner sep=0.75pt]    {$\mathbf{F}( x_{1} ,u_{2})$};
\draw (291.79,182.71) node [anchor=north west][inner sep=0.75pt]    {$f_{1}$};
\draw (58.74,37.95) node [anchor=north west][inner sep=0.75pt]    {$f_{2}$};
\draw (82.63,116.15) node [anchor=north west][inner sep=0.75pt]    {$\mathbf{F}(x_{1} ,u_{1})$};
\draw (113.38,170.97) node [anchor=north west][inner sep=0.75pt]    {$\mathbf{F}(x_{2} ,u_{1})$};
\draw (149.61,149.84) node [anchor=north west][inner sep=0.75pt]  [font=\scriptsize] [align=left]{}; 
\draw (246.72,35.95) node [anchor=north west][inner sep=0.75pt]   [align=left] {nadir};
\draw (96.4,24.6) node [anchor=north west][inner sep=0.75pt]  [color={rgb, 255:red, 208; green, 2; blue, 27 }  ,opacity=1 ]  {$\mathcal{P}^{*}_{\mathbf{m}}( u_{1})$};

\end{tikzpicture}}
    \caption{Illustration of PEHVI. The shaded region in cyan corresponds to the HV improvement corresponding to the choice $x_1$, the fuchsia to $x_2$}
    \label{fig:pehvi_illustration}
\end{figure}
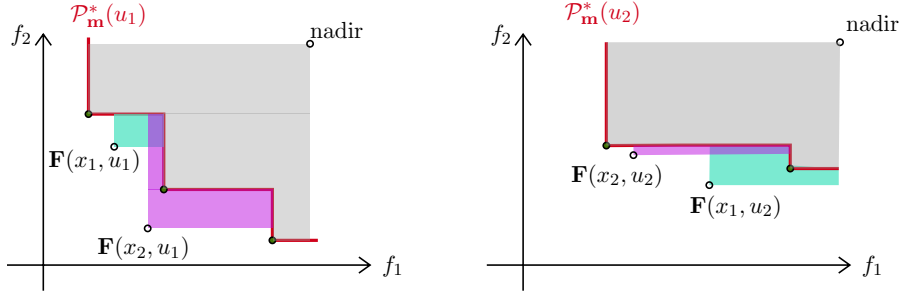

In order to get a better approximation of the Conditional Pareto Fronts and
Conditional Pareto Sets, we can rewrite the EHVI by conditioning on $u \in
\Uspace$, so that we can measure the improvement in the hypervolume of the
dominated region for a given $u$ when adding $\Fb(x, u)$. This conditioned
improvement can be written as
\begin{equation}
    \mathrm{HVI}(\Fb(x, u), \hat{\ParFront}(u)) =  \HV(\hat{\ParFront}(u)\cup \{\Fb(x, u)\}) - \HV(\hat{\ParFront}(u))\, ,
\end{equation}
as illustrated \cref{fig:pehvi_illustration}.

Given that $\Fb$ is a GP, taking the expected value yields the Profile-EHVI,
abbreviated \mbox{PEHVI}, which can be seen as a multiobjective generalization
of the criterion introduced in \cite{ginsbourger_bayesian_2014}:
\begin{align}
        \mathrm{PEHVI}(x, u) =  \Ex_{\mathbf{F}(x,u)}\left[\mathrm{HV}\left(\hat{\ParFront}(u) \cup \Fb(x, u)\right) - \mathrm{HV}(\hat{\ParFront}(u))\right] \, ,
        \label{eq:pehvi}
\end{align}
Optimizing this acquisition function on the joint space $\Xspace \times \Uspace$
gives the next point to evaluate and to add to the design, as explained
\cref{alg:bo_loop}:        
\begin{align}
        (x_{n+1},u_{n+1}) &= \argmax_{(x, u) \in \Xspace \times \Uspace} \mathrm{PEHVI}(x, u) \,.
        \label{eq:pehvi_optim}
\end{align}

\subsection{EHVI under uncertainties}
In order to take into account the random nature of $U$, we can look at the
average of the PEHVI. that we will call ``Integrated EHVI'', to avoid confusion
with the classical EHVI:
\begin{align}
    \mathop{\mathrm{IEHVI}}(x) 
    &= \int_{\Uspace} \Ex_{\Fb(x, u)}\left[\mathrm{HVI}(\Fb(x, u), \hat{\ParFront}(u))\right]p_U(u) \,\mathrm{d}u \label{eq:iehvi_integrals} \, .
\end{align}
This kind of integrated criterion can also be found in noisy extensions of the
EHVI as done in \cite{daulton_parallel_2021} \vicadd{for stochastic simulators}. The integrated criterion can be
used as an acquisition function that we can maximize with respect to the
control variable. Once $x_{n+1}$ has been determined, we sample the next
uncertain variable:
\begin{equation}
\left\{
\begin{array}{rcl}
    x_{n+1} &=& \argmax_{x \in \Xspace} \mathrm{IEHVI}(x) \\
    u_{n+1} &\sim& U
\end{array}\right. \,.
\end{equation}
$x_{n+1}$ is then the point that, when added to the design of experiment, would
increase on average the most the hypervolume of the dominated space averaged with respect to $U$.

In practice, computing the nested integrals of (\ref{eq:iehvi_integrals}) is
untractable analytically. If the inner integral can be computed relatively
efficiently using the same methods of the classical EHVI, the expectation with
respect to $U$ must be approximated using Monte-Carlo and a set of iid samples
of $U$, which can be fixed for all iterations in a Common Random Number fashion
as in \cite{elamri_sampling_2023}, or resampled at every iteration of the BO
loop. 

Nonetheless, for every candidate $x\in\Xspace$, computing the IEHVI still
requires the estimation of the CPF $\hat{\ParFront}(u)$ for each of the samples
of $U$ and then computing the Hypervolume improvement for each of those. This
can get quite expensive if the number of samples is chosen large and if the
estimation of the CPF is expensive.

When comparing the PEHVI with the IEHVI, we traded the optimization of an
integrated criterion with the optimization of a criterion in the joint space
$\Xspace \times \Uspace$ which might be challenging if the dimensionality is too
large, but this criterion benefits from the computationally tractable properties
of the classical EHVI (if we assume that the GP modeling of the different
objective functions are independent), including the availability of the
gradients through analytical derivation \cite{yang_multiobjective_2019} or
automatic differentiation \cite{daulton_parallel_2021}.

In order to take into account both the information on the distribution of $U$,
and to use an acquisition in the joint space $\Xspace\times \Uspace$, we are
also going to compare the heuristic of weighting the PEHVI by $p_U$, the probability
density function of $U$ to retrieve the integrand of
(\ref{eq:iehvi_integrals}):
\begin{equation}
    \mathrm{WPEHVI}(x,u) = \mathrm{PEHVI}(x,u)\cdot p_U(u)\,.
\end{equation}

\subsection{Estimation of the CPF using GP mean}
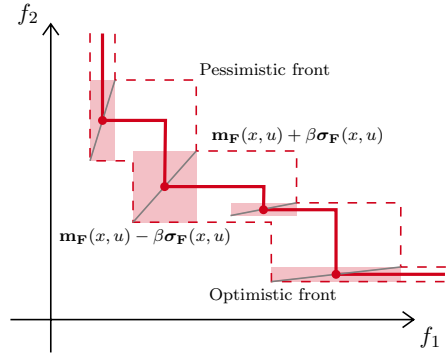
\begin{figure}[ht]
    \centering
    \scalebox{0.8}{\tikzset{every picture/.style={line width=0.75pt}} 

\begin{tikzpicture}[x=0.75pt,y=0.75pt,yscale=-1,xscale=1]

\draw [draw opacity=0][fill={rgb, 255:red, 208; green, 2; blue, 27 }  ,fill opacity=0.25 ]   (426.79,196.42) -- (345.71,205.58) ;
\draw  [draw opacity=0][fill={rgb, 255:red, 208; green, 2; blue, 27 }  ,fill opacity=0.25 ] (345.71,196.42) -- (426.79,196.42) -- (426.79,205.58) -- (345.71,205.58) -- cycle ;

\draw [draw opacity=0][fill={rgb, 255:red, 208; green, 2; blue, 27 }  ,fill opacity=0.25 ]   (361.42,156.17) -- (320.08,164.17) ;
\draw  [draw opacity=0][fill={rgb, 255:red, 208; green, 2; blue, 27 }  ,fill opacity=0.25 ] (320.08,156.17) -- (361.42,156.17) -- (361.42,164.17) -- (320.08,164.17) -- cycle ;

\draw [draw opacity=0][fill={rgb, 255:red, 208; green, 2; blue, 27 }  ,fill opacity=0.25 ]   (298.58,123.67) -- (258.92,168.33) ;
\draw  [draw opacity=0][fill={rgb, 255:red, 208; green, 2; blue, 27 }  ,fill opacity=0.25 ] (258.92,123.67) -- (298.58,123.67) -- (298.58,168.33) -- (258.92,168.33) -- cycle ;

\draw [draw opacity=0][fill={rgb, 255:red, 208; green, 2; blue, 27 }  ,fill opacity=0.25 ]   (247.58,79.17) -- (231.92,129.83) ;
\draw  [draw opacity=0][fill={rgb, 255:red, 208; green, 2; blue, 27 }  ,fill opacity=0.25 ][line width=0.75]  (231.92,79.17) -- (247.58,79.17) -- (247.58,129.83) -- (231.92,129.83) -- cycle ;

\draw  (182,229.4) -- (435,229.4)(207.3,44) -- (207.3,250) (428,224.4) -- (435,229.4) -- (428,234.4) (202.3,51) -- (207.3,44) -- (212.3,51)  ;
\draw [color={rgb, 255:red, 208; green, 2; blue, 27 }  ,draw opacity=1 ][line width=1.5]    (239.75,50) -- (239.75,104.5) -- (279.25,104.5) -- (278.75,146) -- (341,146) -- (340.75,160.5) -- (343.5,160.5) -- (386.25,160.5) -- (386.25,201) -- (454.6,201) ;
\draw [color={rgb, 255:red, 208; green, 2; blue, 27 }  ,draw opacity=1 ] [dash pattern={on 4.5pt off 4.5pt}]  (247.58,50.2) -- (247.58,79.17) -- (298.56,79.17) -- (298.58,123.67) -- (361.56,123.67) -- (361.42,156.17) -- (426.89,156.17) -- (426.79,196.42) -- (460.6,196.42) ;
\draw [color={rgb, 255:red, 208; green, 2; blue, 27 }  ,draw opacity=1 ] [dash pattern={on 4.5pt off 4.5pt}]  (231.92,49.8) -- (231.92,129.83) -- (258.89,129.83) -- (258.92,168.33) -- (345.56,168.33) -- (345.71,205.58) -- (455,205.58) ;
\draw [color={rgb, 255:red, 128; green, 128; blue, 128 }  ,draw opacity=1 ]   (258.92,168.33) -- (298.58,123.67) ;
\draw [color={rgb, 255:red, 128; green, 128; blue, 128 }  ,draw opacity=1 ]   (231.92,129.83) -- (247.58,79.17) ;
\draw [color={rgb, 255:red, 128; green, 128; blue, 128 }  ,draw opacity=1 ]   (320.08,164.17) -- (361.42,156.17) ;
\draw [color={rgb, 255:red, 128; green, 128; blue, 128 }  ,draw opacity=1 ]   (345.71,205.58) -- (426.79,196.42) ;
\draw  [color={rgb, 255:red, 208; green, 2; blue, 27 }  ,draw opacity=1 ][fill={rgb, 255:red, 208; green, 2; blue, 27 }  ,fill opacity=1 ] (237.45,104.5) .. controls (237.45,103.23) and (238.48,102.2) .. (239.75,102.2) .. controls (241.02,102.2) and (242.05,103.23) .. (242.05,104.5) .. controls (242.05,105.77) and (241.02,106.8) .. (239.75,106.8) .. controls (238.48,106.8) and (237.45,105.77) .. (237.45,104.5) -- cycle ;
\draw  [color={rgb, 255:red, 208; green, 2; blue, 27 }  ,draw opacity=1 ][fill={rgb, 255:red, 208; green, 2; blue, 27 }  ,fill opacity=1 ] (276.45,146) .. controls (276.45,144.73) and (277.48,143.7) .. (278.75,143.7) .. controls (280.02,143.7) and (281.05,144.73) .. (281.05,146) .. controls (281.05,147.27) and (280.02,148.3) .. (278.75,148.3) .. controls (277.48,148.3) and (276.45,147.27) .. (276.45,146) -- cycle ;
\draw  [color={rgb, 255:red, 208; green, 2; blue, 27 }  ,draw opacity=1 ][fill={rgb, 255:red, 208; green, 2; blue, 27 }  ,fill opacity=1 ] (338.45,160.5) .. controls (338.45,159.23) and (339.48,158.2) .. (340.75,158.2) .. controls (342.02,158.2) and (343.05,159.23) .. (343.05,160.5) .. controls (343.05,161.77) and (342.02,162.8) .. (340.75,162.8) .. controls (339.48,162.8) and (338.45,161.77) .. (338.45,160.5) -- cycle ;
\draw  [color={rgb, 255:red, 208; green, 2; blue, 27 }  ,draw opacity=1 ][fill={rgb, 255:red, 208; green, 2; blue, 27 }  ,fill opacity=1 ] (383.95,201) .. controls (383.95,199.73) and (384.98,198.7) .. (386.25,198.7) .. controls (387.52,198.7) and (388.55,199.73) .. (388.55,201) .. controls (388.55,202.27) and (387.52,203.3) .. (386.25,203.3) .. controls (384.98,203.3) and (383.95,202.27) .. (383.95,201) -- cycle ;

\draw (306.56,107.4) node [anchor=north west][inner sep=0.75pt]  [font=\scriptsize]  {$\mathbf{m}_{\mathbf{F}}(x,u) +\beta \boldsymbol{\sigma}_{\mathbf{F}}(x,u)$};
\draw (210.96,170.53) node [anchor=north west][inner sep=0.75pt]  [font=\scriptsize]  {$\mathbf{m}_{\mathbf{F}}(x,u) -\beta \boldsymbol{\sigma}_{\mathbf{F}}(x,u)$};
\draw (299.2,67.2) node [anchor=north west][inner sep=0.75pt]  [font=\scriptsize] [align=left] {Pessimistic front};
\draw (305.2,208.4) node [anchor=north west][inner sep=0.75pt]  [font=\scriptsize] [align=left] {Optimistic front};
\draw (437.2,232.6) node [anchor=north west][inner sep=0.75pt]    {$f_{1}$};
\draw (184.8,29.8) node [anchor=north west][inner sep=0.75pt]    {$f_{2}$};

\end{tikzpicture}}
    \caption{Different estimations of the Pareto front using bounding-boxes. The shaded regions corresponds to the bounding boxes}
    \label{fig:beta_front}
\end{figure}

The estimation of the CPF $\hat{\ParFront}(u)$ \vicadd{which is used in the acquisition functions in} (\ref{eq:iehvi_integrals} \ref{eq:pehvi}) can be done in different ways. We propose
here to use a \emph{conservative} approximation based on the GP mean $\mb_\Fb$
of the objective function, as can be found in bounding-box approaches
\cite{rivier_surrogateassisted_2022}, evaluated on a discrete set of input
points $\mathcal{X}_{\text{pareto}} \subset \Xspace$
\begin{align}
    \hat{\ParFront}(u) = \text{non-dominated points of }\left\{ \mb_{\Fb}(x, u) + \beta \boldsymbol{\sigma}_{\Fb}(x, u) \mid x \in \mathcal{X}_{\text{pareto}}\right\}\,, \label{eq:parfront_beta}
\end{align} as illustrated \cref{fig:beta_front}.
\vicadd{Since the different objectives are modelled independently, the covariance matrix of the predictions is diagonal, and we define the vector $\boldsymbol{\sigma}_{\Fb}$ as its diagonal, and represents thus the prediction variance of each objective independently.}

This estimation relies on two additional parameters:
\begin{itemize}
    \item The value of $\beta\in\mathbb{R}$. If $\beta>0$, we have a
    ``pessimistic" front, since we are being conservative with our estimate of
    the CPF. Conversely, if $\beta <0$, then we have a ``optimistic" estimation
    of the CPF.
    \item $|\mathcal{X}_{\mathrm{pareto}}|$, i.e. the number of
    points used to approximate the Pareto front
\end{itemize}
Their influence are studied more closely in the next section.
\vicadd{When computing the probability of coverage, we can estimate the CPF using (\ref{eq:parfront_beta}) and $\beta=0$}

\section{Numerical Experiments}
We are going to compare the three different acquisition function introduced
above, namely the PEHVI, the WPEHVI and the IEHVI. All those methods are to be
compared with a random filling of the input space. 

In order to compare those, we are going to introduce the appropriate metrics
first. In this work, we will estimate the CPS and CPF using a plug-in approach,
by using the mean if the Gaussian process conditioned on the design of
experiments as a replacement of the true function. For a given $u\in\Uspace$,
and a set $\mathcal{X}_{\text{test}}\subset\Xspace$ \vicadd{of testing points different from $\Xspace_{\text{pareto}}$}, we can compare the CPF $\ParFront(u)$
computed as the non-dominated points of $\{\fb(x, u) \mid
x\in\mathcal{X}_{\text{test}}\}$ and its plug-in estimation $\hat{\ParFront}(u)$
which is the set of the non-dominated points of $\{\mb_{\Fb}(x, u) \mid
x\in\mathcal{X}_{\text{test}}\}$ \vicadd{, which follows (\ref{eq:parfront_beta}) computed with $\beta=0$}.

\subsection{Metrics used for comparison}

\paragraph{Average Hausdorff distance to the true Pareto front}
We are first going to introduce a metric to compare an estimation of the
Pareto front with its true value. In the objective space, the distance
of a point $\yb$ to a set $\mathcal{A}\subset \mathbb{R}^d$ is defined as
\begin{align}
    d(\yb, \mathcal{A}) = \inf_{\mathbf{z}\in \mathcal{A}} \|\yb - \mathbf{z}\| \,.
\end{align}
In \cite{schutze_using_2012}, the Generational Distance (GD) and Inverse
Generational Distance (IGD) between an estimation of the Pareto front
$\hat{\ParFront}$ and the true Pareto front $\ParFront$ are defined as 
\begin{align}
    \mathop{\mathrm{GD}_p}(\hat{\ParFront}, \ParFront) = \left(\frac{1}{|\hat{\ParFront}|} \sum_{\hat{\yb} \in \hat{\ParFront}} d(\hat{\yb}, \ParFront)^p\right)^{1/p} \, ,
\end{align}
and
\begin{align}
    \mathop{\mathrm{IGD}_p}(\hat{\ParFront}, \ParFront) = \left(\frac{1}{|\ParFront|} \sum_{\yb \in {\ParFront}} d(\yb, \hat{\ParFront})^p\right)^{1/p}\,.
\end{align}
$\mathop{\mathrm{GD}_p}$ is proportional to the average $L_p$ distance of a
point of the approximation of the front to the reference front, while the
$\mathop{\mathrm{IGD}_p}$ is proportional to the average $L_p$ distance of a
point of the reference front to the approximated one. In
\cite{schutze_using_2012}, the authors propose \vicadd{to combine both by defining} the averaged Hausdorff distance,
denoted as $\Delta_p$ as a quality metric with respect to the true Pareto front:
\begin{align}
    \Delta_p(\hat{\ParFront}, \ParFront) = \max(\mathop{\mathrm{GD}_p}(\hat{\ParFront}, \ParFront), \mathop{\mathrm{IGD}_p}(\hat{\ParFront}, \ParFront))
\end{align}

This averaged Hausdorff distance can be used for a given $u$ to compare the true
CPF with its plug-in estimation: $\Delta(u) = \Delta_2(\hat{\ParFront}(u),
\ParFront(u))$. We can then compare the distributions of $\Delta(U)$, estimated using a finite set of samples $\{u_i\}_{1\leq i\leq n_u}$ depending on the design of experiment obtained using the different acquisition functions. \vicadd{The distribution of the average Hausdorff distance is represented \cref{fig:hausdorff_5x5_unif,fig:hausdorff_5x5_gaussian}}

\paragraph{Comparison on the probability of coverage}
In a similar way, we can estimate the probability of coverage using a set $\{u_i\}_{1\leq i\leq n_u}$ of $n_u$ samples of $U$, defined for $x\in\mathcal{X}_{\mathrm{test}}$ as
\begin{align}
   p_\mathrm{cov}(x) = \frac{1}{n_u}\sum_{i=1}^{n_u} \mathbbm{1}_{\ParSet(u_i)}(x), \qquad
      \hat{p}_\mathrm{cov}(x) = \frac{1}{n_u}\sum_{i=1}^{n_u} \mathbbm{1}_{\hat{\ParSet}(u_i)}(x)
\end{align}
We can then compute the average $L_2$ distance over the whole input space: 
\begin{equation}
    \ell(p_\mathrm{cov},  \hat{p}_\mathrm{cov}) = \frac{1}{|\mathcal{X}_{\mathrm{test}}|}\sum_{i=1}^{|\mathcal{X}_{\mathrm{test}}|} (p_\mathrm{cov}(x_i)-  \hat{p}_\mathrm{cov}(x_i))^2
\end{equation}
\subsection{Analytical toy problems}
\label{sec:4d_prob}
\subsubsection{4D problem}
\paragraph{Experimental setting}
We are going to see the influence of the criteria introduced before on the
multiobjective optimization problem of toy functions.
We are first interested in a toy problem where both the control variable and the
uncertain variable are two-dimensional:
\begin{equation}
    \Xspace = [0, 1] \times [1, 2], \quad  \Uspace = [2, 3] \times [3, 4]\,.
    \end{equation}
The objective function is defined analytically as
\begin{equation}
        \begin{array}{rcl}
         \fb_{2\times 2}: \Xspace \times \Uspace& \longrightarrow& \mathbb{R}^2 \\
         (x,u) & \longmapsto & \begin{pmatrix}
             (x_1 - u_1 + 2)^2 + (x_2 - u_2 + 2)^2 + 5 u_1\\
             (x_1 - x_2 + 1)^2 + (x_1x_2 - u_1 + 1.5)^2 + 5u_2
         \end{pmatrix}\,.
    \end{array}
    \label{eq:def_f2x2}
\end{equation}

\begin{figure}
    \centering
    \includegraphics[width=.75\linewidth]{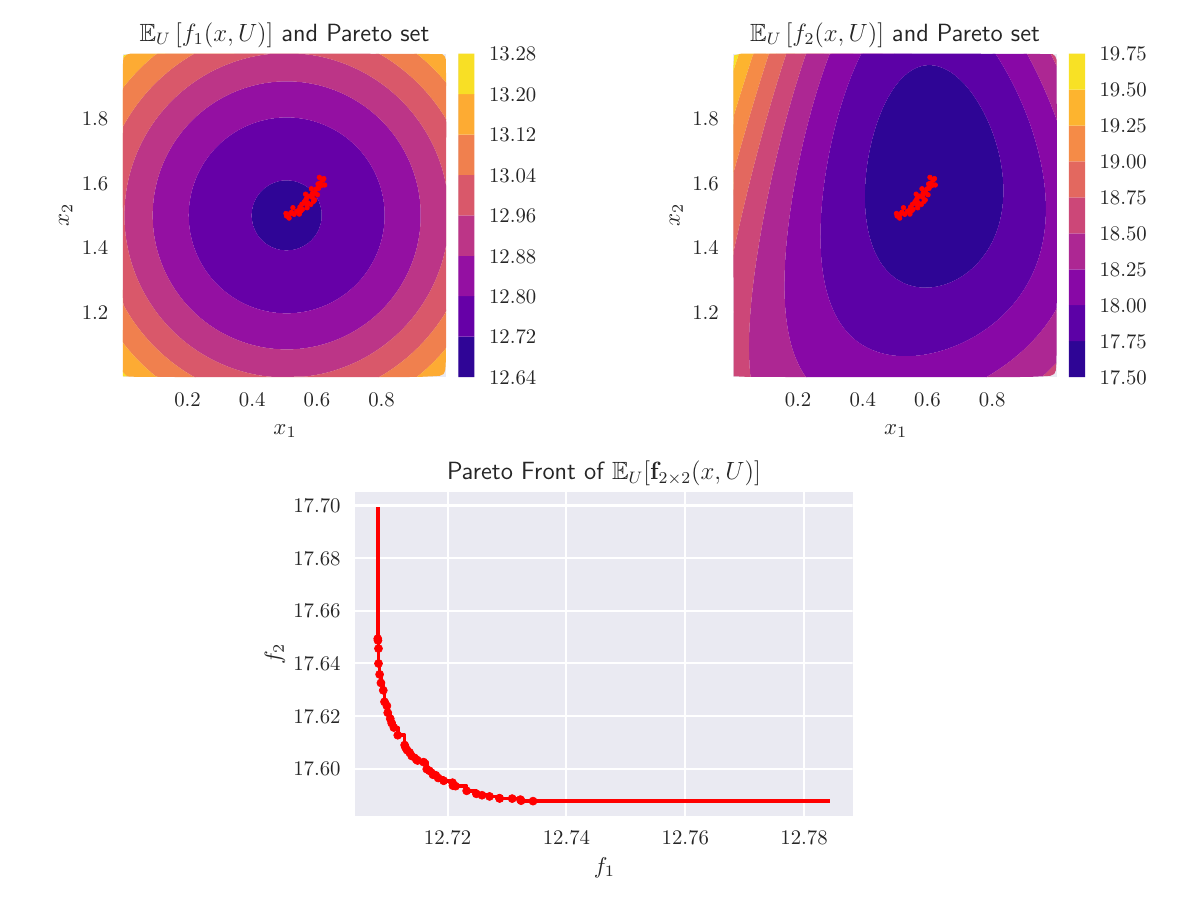}
    \caption{Mean objective multiobjective optimization: the mean of the first and second objective are shown on the \vicreplace{leftmost and middle figures}{top figures}. The Pareto front is shown on the \vicreplace{rightmost}{bottom} figure.}
    \label{fig:mean_obj}
\end{figure}

If we assume that $U \sim \mathrm{Unif}(\Uspace)$, the mean objective is
\begin{equation}
    \Ex_U[\fb_{2\times 2}(x, U)] = \begin{pmatrix}
        x_1(x_1 - 1) + x_2(x_2 - 3) + \frac{91}{6} \\
        x_1^2 x_2^2 + x_1^2 - 4 x_1x_2 + 2x_1 + x_2^2 - 2x_2 + 235/12
    \end{pmatrix}\,.
\end{equation}
Due to the simplicity of the problem, we can visualize some quantities easily:
the averaged MO problem is illustrated \cref{fig:mean_obj}, while
\cref{fig:prob_cpf} shows the estimated probability of coverage, as defined
(\ref{eq:prob_coverage}), computing numerically using a large number of 
evaluations of the true function.


\Cref{fig:prob_cpf} shows a behavior which is quite different from the Pareto
front of the objective means. Indeed, we can see that the regions of interest
differ: regions of relatively high probability in \cref{fig:prob_cpf} are
located around the diagonal, and especially for $x_1$ close to $0.4$. When
considering the mean of the objectives, we can see that the Pareto set is located
mostly around $x_1\approx 0.5$.

\begin{figure}[ht]
\centering
\begin{subfigure}{0.7\textwidth}
    \centering
    \includegraphics[width=.8\textwidth]{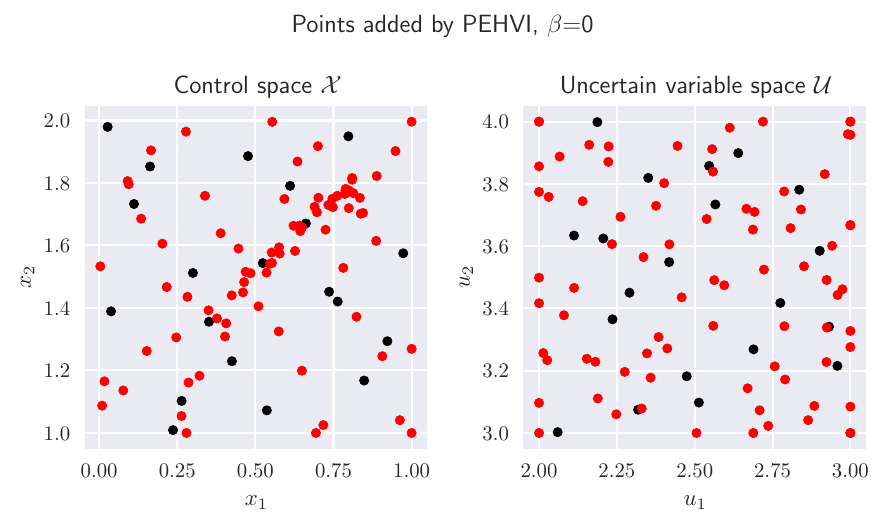}
    \caption{Design of experiment obtained using the PEHVI with $\beta=0$}
    \label{fig:pehvi_beta_00}
\end{subfigure}
\begin{subfigure}{0.7\textwidth}
    \centering
    \includegraphics[width=.8\textwidth]{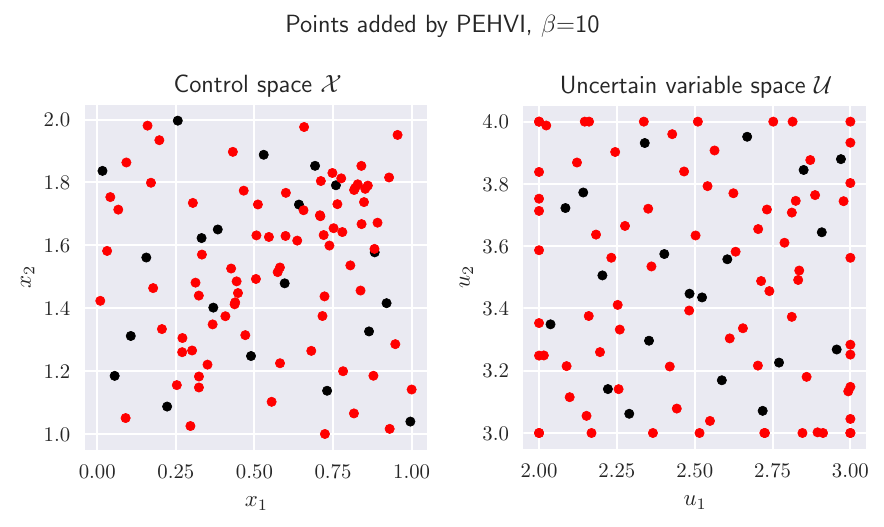}
    \caption{Design of experiment obtained using the PEHVI with $\beta=10$}
    \label{fig:pehvi_beta_10}
\end{subfigure}
\caption{Final designs of experiment using PEHVI, for different values of
$\beta$ for the estimation of $\ParFront(u)$ \vicadd{ as described (\ref{eq:parfront_beta}}). The initial design of
experiment is represented in black, the added points are in red}
\label{fig:pehvi_beta}
\end{figure}

\paragraph{Influence of the hyperparameters}
In a parallel with the classical EI \cite{jones_efficient_1998}, the $\beta$ parameter controls the
incumbent, i.e. the ``best" value so far, that is the current estimation of the
Pareto front. Numerical experiments suggest that taking a \vicreplace{moderate number of}{few hundreds}
points for the estimation, \vicreplace{and a}{coupled with a} ``pessimistic'' \vicadd{($\beta >0$)} estimation of the front  leads to
better results. Indeed, such choices push toward more exploration of the input
space, as seen on \cref{fig:pehvi_beta} which shows the points added using the
PEHVI acquisition function.

\vicadd{Considering the $x$ components of those added points (left plots), they are mostly added along the diagonal, which is a region of interest, as seen on \cref{fig:prob_cpf}. The $u$ components of the added points, obtained through the optimization of the criterion, are spread more uniformly within the input space.}
However, it seems that this acquisition intensifies well, but does not explore
sufficiently the input space when $\beta$ is zero. A better exploration around
this diagonal is obtained using a larger value for $\beta$.

\begin{figure}[ht]
    \centering
    \includegraphics[scale=0.8]{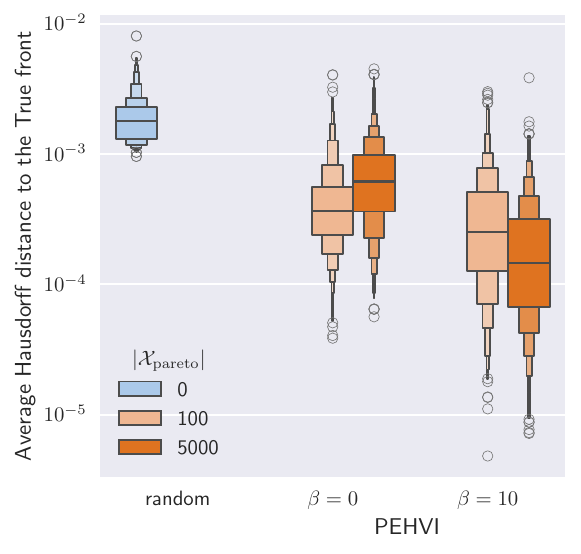}
    \caption{Averaged Hausdorff distance comparing the estimated Pareto front
    using the GP mean, depending on the hyperparameter $\beta$ and the number of
    points in $\mathcal{X}_{\text{pareto}}$ \vicadd{used during the optimization of the PEHVI for the 4D problem}. Performance of the random design
    for reference \vicadd{(corresponding to  $|\mathcal{X}_{\text{pareto}}|=0$) in the figure.}}
    \label{fig:boxplot_hausdorff}
\end{figure}

This is also illustrated \cref{fig:boxplot_hausdorff}, where we can see the
difference in the averaged Hausdorff distance for different values of $\beta$,
and a different number of points for the estimation of the CPF. Since the PEHVI
measures an improvement in  the HV compared to an estimated CPF, having a
relatively conservative estimation through a positive value of $\beta$ or a small
$|\mathcal{X}_{\text{pareto}}|$ helps improve the performance of the algorithm.

\subsubsection{10D problem}
\paragraph{Experimental setting}
We will consider now the following multiobjective optimization problem where the control space
$\Xspace$ has dimension 5, while the uncertain space is also 5-dimensional:
\begin{equation}
    \Xspace = [0, 1]^5, \quad  \Uspace =  [0, 1]^5\,,
    \end{equation}
and the objective function is defined analytically as
\begin{equation}
        \begin{array}{rcl}
         \fb_{5\times 5}: \Xspace \times \Uspace& \longrightarrow& \mathbb{R}^2 \\
         (x,u) & \longmapsto & \begin{pmatrix}
             (\sum_{i=1}^5 x_i + u_1 + u_2 + u_3  - u_4 + u_5 - 5)^2 \\
             (\sum_{i=1}^5 x_i + u_1 + u_2 + u_3 + u_4 - u_5 - 5)^2
         \end{pmatrix}\,.
    \end{array}
\end{equation}

We consider two different modeling for the uncertain variable leading to two different problems of multiobjective optimization under uncertainties:
\begin{itemize}
    \item Problem 10d with $U\sim\mathrm{Unif}(\mathcal{U})$, and
    \item Problem 10d bis with $U\sim\mathcal{N}(u_0, \Sigma_0)$
\end{itemize}
where $\Sigma_0$ is diagonal, with $10^{-1}$ on the diagonal, and $u_0$ is the
center of the domain. 

On both problems, we sample an initial design in the joint
space of 100 points \vicadd{according to the usual rule of thumb $10 \cdot \dim(\Xspace\times\Uspace)$}, and add 400 points sequentially according to the criteria
introduced above. \vicadd{For each of these methods, 10 replications were performed in order to account partially for the randomness in the methods.}

\begin{figure}
    \centering
    \begin{subfigure}{0.49\textwidth}
    \includegraphics[width=.95\linewidth]{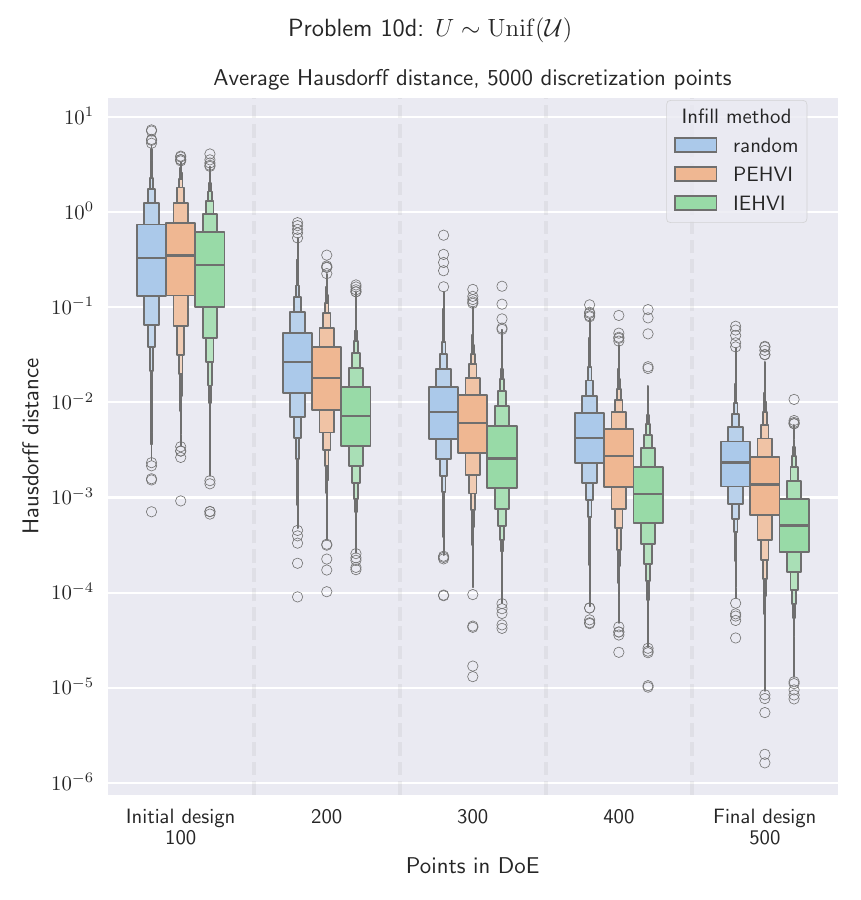}
    \caption{\vicadd{10d problem}}
    \label{fig:hausdorff_5x5_unif}
    \end{subfigure}
    \begin{subfigure}{0.49\textwidth}
            \centering
    \includegraphics[width=.95\linewidth]{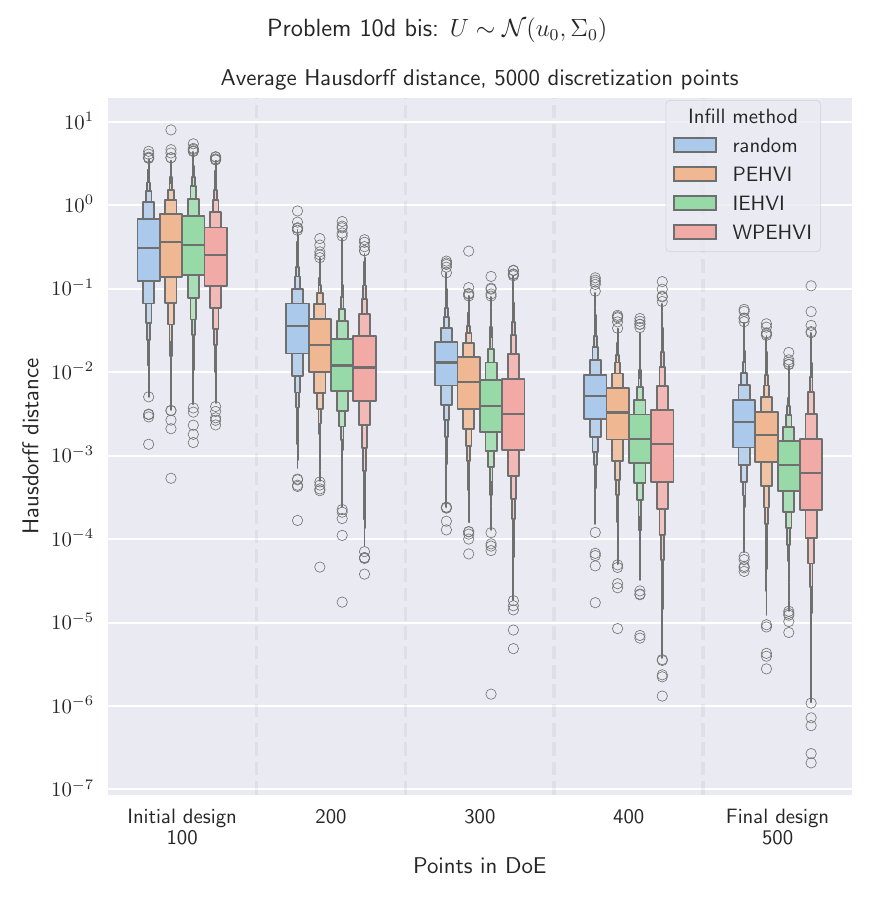}
    \caption{10d problem bis}
    \label{fig:hausdorff_5x5_gaussian}
    \end{subfigure}
    \caption{Average Hausdorff distance for 512 samples of $U$, and 10 replications for each experiment, and  evaluate the different CPF}
    \label{fig:avg_hausdorff}
\end{figure}

\begin{figure}
    \centering
    \begin{subfigure}{0.49\textwidth}
    \includegraphics[width=.95\linewidth]{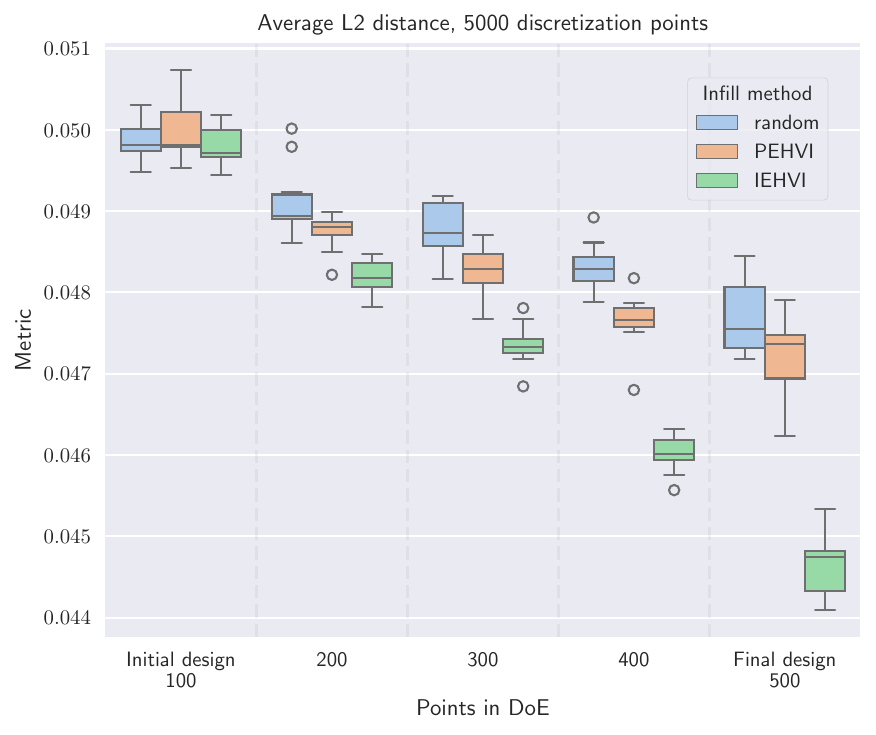}
    \caption{\vicadd{10d problem}}
    \label{fig:l2_5x5_unif}
    \end{subfigure}
    \begin{subfigure}{0.49\textwidth}
            \centering
    \includegraphics[width=.95\linewidth]{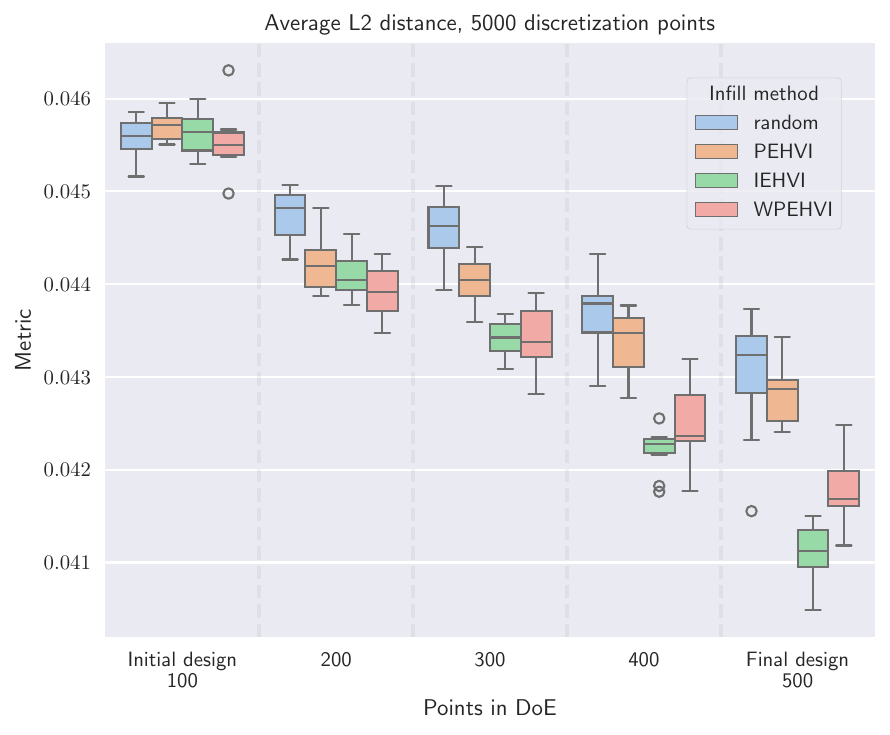}
    \caption{10d problem bis}
    \label{fig:l2_5x5_gaussian}
    \end{subfigure}
    \caption{Averaged $L_2$ distance (over 512 samples of $U$) for each problem. 10 replications for each method and $|\Xspace_{\mathrm{test}}| = 5000$}
    \label{fig:avg_l2}
\end{figure}

\paragraph{Results}
\Cref{fig:avg_hausdorff} show the averaged Hausdorff distance of the plug-in
estimates depending on the acquisition function used, while \cref{fig:avg_l2} shows the averaged $L_2$ distance between the estimated probability of coverage. For the problem with
uniformly distributed uncertainties, the PEHVI performs better than the random
design, while the IEHVI performs better than the PEHVI. The same ordering
appears for the non-uniform case, as shown on \cref{fig:hausdorff_5x5_gaussian}\vicadd{, even though the PEHVI is distribution-agnostic. This can be explained by the space-filling properties of the acquisition function, as illustrated \cref{fig:pehvi_beta}. }When comparing with the WPEHVI, the estimated
distribution of the averaged Hausdorff distance shows a larger variance, with the
median close to the one of the IEHVI, but more extreme values skew the
distribution, leading to a mean close to the one of the random design.

In both problems, the IEHVI gives significantly
better results than the other methods. \vicdelete{According to this metric, for the problem with
uniformly distributed $U$, the difference between PEHVI and random design is
less significant, but for the normally distributed $U$, the difference between
those two methods are not significant, but all the other differences are.}\vicadd{Comparing the PEHVI and the random design, we can see that for Problem 10d (with uniformly distributed $U$), PEHVI still seems to give better results than a random design. For Problem 10d bis, PEHVI and random design give similar performances, while IEHVI and WPEHVI perform better.}


Based on these results, the IEHVI seems to provide consistently better results
than the other acquisition functions. However, each of its evaluation requires
multiple estimations of CPF, for each of the samples of $U$ chosen to evaluate
the integral. This may become tedious as the dimension of the problem increases,
and the PEHVI or WPEHVI can be simpler alternatives to implement.
\clearpage
\subsection{Design of a cabin using EnergyPlus}
As  an illustration, we are going to apply this principle of multiobjective
optimization to a toy problem of cabin design. The cabin is a simple
construction of 9 \unit{\meter\squared}, with a window, a door, a heating unit
and a cooling unit. The walls are made of a concrete layer, and an insulator
layer, while the floor is made of concrete. We chose to model the control
parameters as the thicknesses of the different layers, while the uncertain
parameters represent the individual user preferences regarding the temperature
setpoints, or their tendency to ventilate the cabin. Those parameters are
described \cref{tab:cabin_design}.

\begin{table}[ht]
    \centering
    \begin{tabular}{lccccl}\toprule
                  &Physical parameter & Unit  & Space\\\midrule
    \multirow{3}{*}{Control: $x$} &  Thickness of concrete wall   & \unit{\meter}   & $[0.05, 0.30]$    \\
        &Thickness of concrete floor &  \unit{\meter}  & $[0.05, 0.30]$    \\
        &Thickness of wall insulator  & \unit{\meter}  & $[0.05, 0.30]$   \\\midrule
        \multirow{3}{*}{Uncertain: $u$} & Air infiltration   & {Air change/hour} & $\text{Unif}\left([1.0, 4.0]\right)$   \\
        & Temperature setpoint for heating &  \unit{\degreeCelsius} & $\text{Unif}\left([18.0, 22.0]\right)$  \\
        & Temperature setpoint for cooling&   \unit{\degreeCelsius}  & $\text{Unif}\left([24.0, 28.0]\right)$ \\\bottomrule
    \end{tabular}
    \caption{Configuration of the cabin design problem}
    \label{tab:cabin_design}
\end{table}
Using
EnergyPlus\footnote{\href{https://energyplus.net/}{https://energyplus.net/}}, we are able to compute three different quantities that are to be
optimized:
\begin{itemize}
    \item the total energy needed for the heating and cooling
    units,
    \item the comfort index for the occupants of the cabin computed using
    \texttt{pythermalcomfort} \cite{tartarini_pythermalcomfort_2020},
    \item the cost of the materials (proportional to the volume of each material
    needed).
\end{itemize}

Based on this, we consider two different Multiobjective Optimization problems,
with three (\ref{eq:oreni_3}) or two objectives (\ref{eq:oreni_2}):
\begin{align}
(x,u)&\longmapsto  \left(\text{Energy}(x,u), -\text{Comfort}(x,u), \text{Cost}(x,u)\right)\,,\text{ and} \label{eq:oreni_3}\\
(x,u) &\longmapsto \left(\text{Energy}(x,u), -\text{Comfort}(x,u)\right) \,.\label{eq:oreni_2}
\end{align}

Based on an initial design of 60 points in $\Xspace\times\Uspace$, we
constructed the initial GP regression models, and applied the procedure
explained above for 140 iterations using the IEHVI, in order to reach a total number of
simulations of 200. \vicadd{We choose to model each objective independently, with a Matérn 5/2 kernel using GPytorch \cite{gardner_gpytorch_2018}. We used Botorch \cite{balandat_botorch_2020} to implement the different acquisition functions.}

\begin{figure}[ht]
    \centering
    \begin{subfigure}{0.49\textwidth}
        \includegraphics[width=\textwidth]{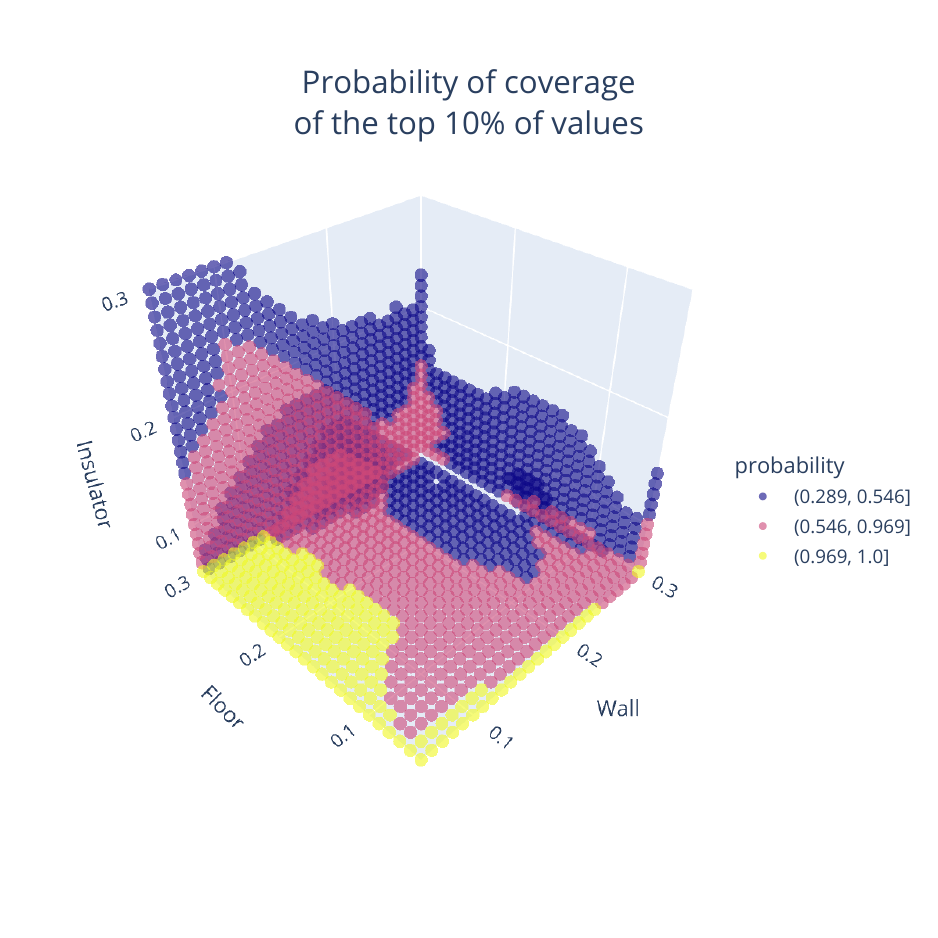}
        \caption{All objectives are considered (\ref{eq:oreni_3})}
        \label{fig:oreni3}
    \end{subfigure}
    \begin{subfigure}{0.49\textwidth}
        \includegraphics[width=\textwidth]{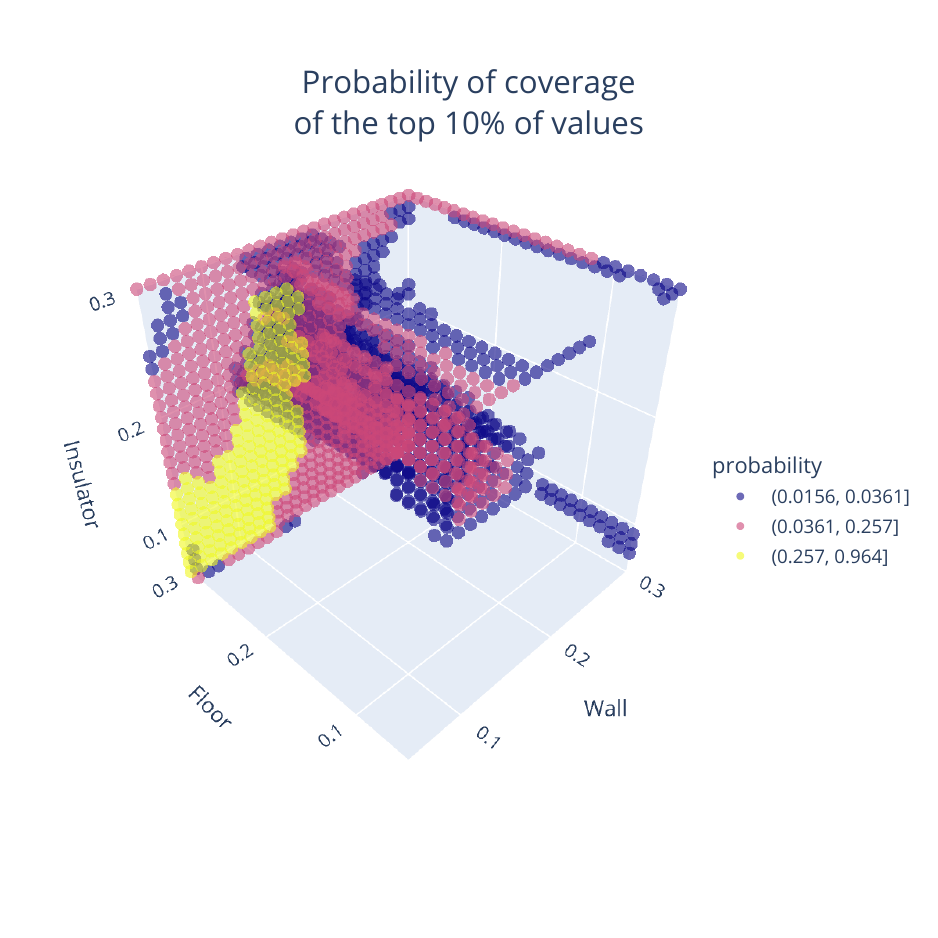}
        \caption{Only the energy spent (for AC and heating units) and the comfort index are considered (\ref{eq:oreni_2})}
        \label{fig:oreni2}
    \end{subfigure}
    \caption{Estimated probability of coverage using the metamodel constructed using the PEHVI acquisition function. \vicadd{The axes of the hypercube are the values of the control variable $x$. The yellow dots are the designs corresponding to the top 1\% of values, the red dots to the top 5\%, and the blue to the top 10\%.}}
    \label{fig:oreni}
    \end{figure}

The two plots of \cref{fig:oreni} show the estimated
probability of coverage, computed using the GP prediction: the probability has
been estimated using 1024 samples of $U$, where the \vicreplace{condidates}{candidate designs} are
arranged on a regular grid of $26^3=17576$ points. \vicadd{For clarity, we only displayed the points with a probability of coverage in the top 10\% of values}.
\vicdelete{The upper one represents the situation where all objectives are considered
simultaneously, while the bottom one shows the probability of coverage when
considering only the comfort index and the total energy needed.}
\vicadd{\Cref{fig:oreni3} represents the situation where all objectives are considered
simultaneously, while \cref{fig:oreni2} one shows the probability of coverage when
considering only the comfort index and the total energy needed.}

\begin{figure}[ht]
    \centering
    \includegraphics[width=0.6\textwidth]{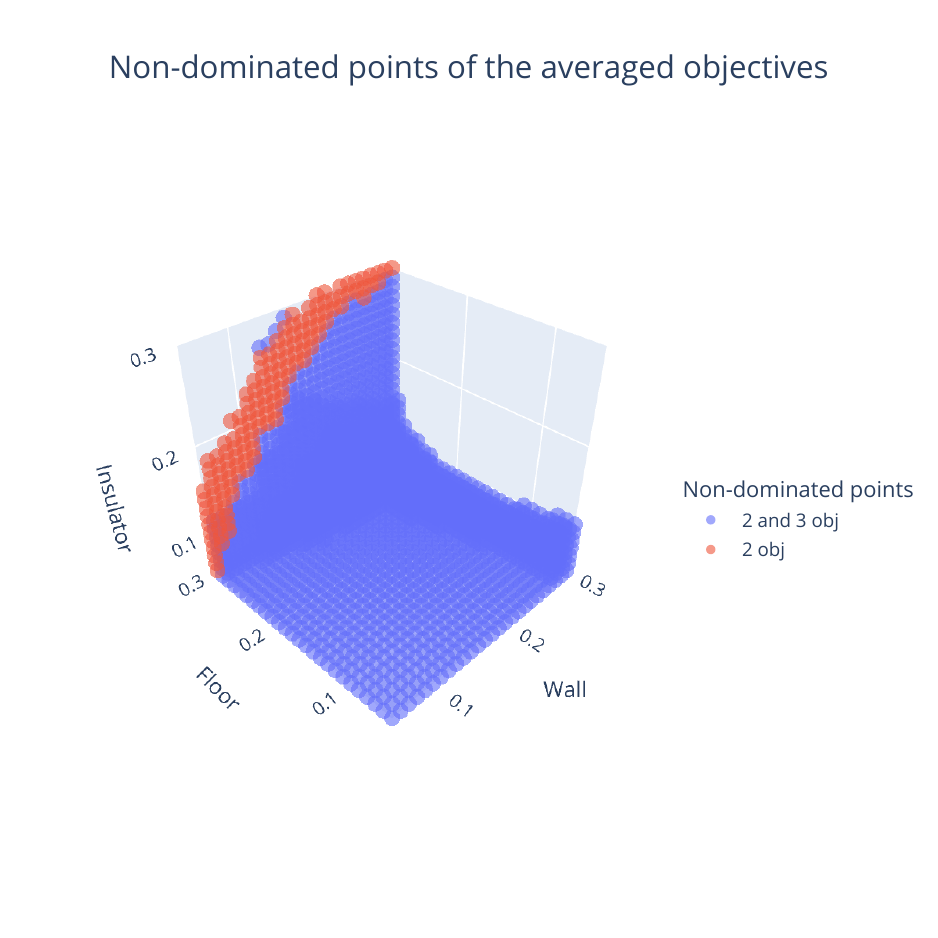}
    \caption{Non-dominated points of the mean of the objectives. The red points correspond to non-dominated points when considering two objectives, the blue ones when considering all three objectives}
    \label{fig:mean_obj_oreni}
\end{figure}

For comparison purposes, we also constructed also a GP model based on a design
 of 1024 points in order to compute the Pareto Front of
the mean of the objectives, as in the multiobjective optimization problem
described (\ref{eq:mean_obj_pareto}). The non-dominated points \vicadd{of the expected value of the objectives, as in (\ref{eq:mean_obj_pareto}}) are
represented~\cref{fig:mean_obj_oreni}, for the problems with two and three
objectives. \vicadd{The points in blue, which correspond to the Pareto set of the expected value of (\ref{eq:oreni_3}), are located mostly on the boundary of the domain, where the thickness of the insulator layer is at its minimum (bottom face of the cube), and where the thickness of the floor is at its maximum (back right face of the cube). The point in red, corresponding to the Pareto set of the expected value of (\ref{eq:oreni_2}) are where the floor thickness is at its maximum, and and when the thickness of the insulator layer and the concrete layer in the wall increase simultaneously.}

When considering three objectives, most non-dominated points are located on the
plane corresponding to a maximal floor thickness, or on the place corresponding
to a minimal insulator thickness. When removing the cost from considerations by
taking only two objectives, all the non-dominated points correspond to a maximal
floor thickness, while a tradeoff appear between insulator thickness and wall
thickness.

We can see that both methods, i.e. looking at the mean of the objectives, and the
probability of coverage lead to similar conclusions on the design of the cabin,
in the sense that the same regions of interest seem to be identified. One
advantage in that case of the probability of coverage is that it highlights the
different behavior under uncertainties of the solutions brought by optimization
of the averaged objectives.

\section{Conclusion and Perspectives}
In this work, we introduced the notions of Conditional Pareto Fronts and
Conditional Pareto Sets which can be used in order to get insight on the
distribution of solutions of the problem of multiobjective optimization. We
propose to use the coverage probability of the CPS in order to sort the
potential designs considered. Since this probability is expensive to compute, we
can use surrogate models based on Gaussian Processes and an Active Learning
approach to improve its ability to predict the different quantities of interest.
\vicadd{Some further work could be investigated in order to improve the acquisition functions:}
\vicreplace{The}{the} IEHVI relies on sampling to select the uncertain variable. \vicadd{In order to tie this selection more tightly to the procedure, one}\vicdelete{One} could derive
a method based on Stepwise Uncertainty Reduction, as done in
\cite{elamri_sampling_2023} in order to select \vicreplace{the uncertain variable.}{it.}

Instead of just considering the probability with which a point belongs
or not to the CPS, one could look at introducing a measure of distance to the
front, using for instance non-dominated sorting as in NSGA-II
\cite{deb_fast_2002}, or by using measures of quality of Pareto front
approximation \cite{audet_performance_2021}, in order to get directly a set of
points which could be considered as a robust counterpart of the Pareto front.

On a more general note, the acquisition functions introduced in this work rely
on the ability to model the objective functions in the joint space
$\Xspace\times \Uspace$, and on the assumption that we know the distribution of
the uncertain variable. In that sense, this work could be extended to the case
where we do not have access directly to the distribution of $U$, but only have
access to a limited number of samples of $U$, in an unknown space $\Uspace$, and
thus need to fit a GP model for every available sample of $U$. Removing even
more assumptions on the distribution of $U$ could lead to approaching the problem
from a distributionally robust optimization point of view
(see \cite{lin_distributionally_2022}) by introducing an ambiguity set on the
distribution of $U$. Even though the PEHVI acquisition function shows mixed
results in the experiments introduced above, the PEHVI can be derived also from
the IEHVI by considering the maximization of the IEHVI with an ambiguity set on
the distribution of $U$.

\section*{Acknowledgements}
\vicadd{We thank the anonymous reviewers, whose comments improved the paper.}
This work has benefited from the expertise of people of the Tipee
platform in order to apply our methods
to the cabin design problem using EnergyPlus. This work was granted access to the HPC resources
of PMCS2I (Pôle de Modélisation et de Calcul en Sciences de l’Ingénieur de
l’Information) of École Centrale de Lyon, Écully, France. 
\bibliographystyle{apalike}
\bibliography{MOOUU.bib}
\end{document}